\documentclass{amsart}
\usepackage{amssymb,amscd}

\numberwithin{equation}{section}

\newtheorem{prop}{Proposition}
\newtheorem{theorem}[prop]{Theorem}

\newtheorem{lemma}[prop]{Lemma}
\newtheorem{conjecture}[prop]{Conjecture}

\theoremstyle{definition}
\newtheorem{definition}[prop]{Definition}
\newtheorem{example}[prop]{Example}
\newtheorem{remark}[prop]{Remark}

\numberwithin{prop}{section}

\newfont{\germ}{eufm10}
\newfont{\slsmall}{cmsl8}
\newfont{\bfsl}{cmbxsl10}

\newcommand\Bd{B^\dagger}
\newcommand\Bmin{B_{\min}}
\newcommand\bb{\mbox{\bfsl b}}
\newcommand\bd{b^\dagger}
\newcommand\bt{\tilde{b}}
\newcommand\Cfin{{\mathcal C}^{\mathrm{fin}}}

\newcommand\et[1]{\tilde{e}_{#1}}
\newcommand\ft[1]{\tilde{f}_{#1}}
\newcommand\geh{\mbox{\germ g}}
\newcommand\La{\Lambda}
\newcommand\la{\lambda}
\newcommand\lev{\mbox{\sl lev}\,}
\newcommand{\na}{\natural}
\newcommand\ol{\overline}
\newcommand\ot{\otimes}
\renewcommand\P{{\mathcal P}}
\newcommand\PB{{\mathcal P}(\pb,B)}
\newcommand\pb{\mbox{\bfsl p}}
\newcommand\Pcl{P_{cl}}
\newcommand\Pcll{(\Pcl^+)_l}
\newcommand\simarrow{\stackrel{\sim}{\rightarrow}}
\newcommand\Uq{U_q(\geh)}
\newcommand\Uqp{U'_q(\geh)}
\newcommand\veps{\varepsilon}
\newcommand\vphi{\varphi}
\newcommand\wt{\mbox{\sl wt}\,}
\newcommand\wtb{\mbox{\bfsl wt}\,}
\newcommand\wts{\mbox{\slsmall wt}\,}
\newcommand\Z{\mathbb{Z}}
\newcommand\Zn{\Z_{\ge0}}

\begin{document}

\title{A Tensor Product Theorem Related to \\ Perfect Crystals}

\author[M.~Okado]{Masato Okado}
\address{Department of Informatics and Mathematical Science,
Graduate School of Engineering Science, Osaka University,
Toyonaka, Osaka 560-8531, Japan}
\email{okado@sigmath.es.osaka-u.ac.jp}

\author[A.~Schilling]{Anne Schilling}
\address{Department of Mathematics, University of California, One Shields
Avenue, Davis, CA 95616-8633, U.S.A.}
\email{anne@math.ucdavis.edu}

\author[M.~Shimozono]{Mark Shimozono}
\address{Department of Mathematics, 460 McBryde Hall, Virginia Tech,
Blacksburg, VA 24061-0123, U.S.A}
\email{mshimo@math.vt.edu}


\begin{abstract}
Kang et al. provided a path realization of the crystal graph of a
highest weight module over a quantum affine algebra, as certain
semi-infinite tensor products of a single perfect crystal. In this
paper, this result is generalized to give a realization of the
tensor product of several highest weight modules. The underlying
building blocks of the paths are finite tensor products of several
perfect crystals. The motivation for this work is an
interpretation of fermionic formulas, which arise from the
combinatorics of Bethe Ansatz studies of solvable lattice models,
as branching functions of affine Lie algebras. It is shown that
the conditions for the tensor product theorem are satisfied for
coherent families of crystals previously studied by Kang,
Kashiwara and Misra, and the coherent family of crystals
$\{B^{k,l}\}_{l\ge 1}$ of type $A_n^{(1)}$.
\end{abstract}

\maketitle


\section{Introduction}
In the seminal paper \cite{KMN1}, the crystal graph $B(\la)$ of
the irreducible integrable module of highest weight $\la$ over a
quantum affine algebra $\Uq$, is realized as a subset
$\P(\pb^{(\la)},B)$ of the semi-infinite tensor product of copies of
the crystal graph $B$ of a single finite-dimensional
$\Uqp$-module, where $\Uqp$ is the subalgebra of $\Uq$
corresponding to the derived subalgebra $\geh'$ of $\geh$. The
suitable finite crystals $B$ are called perfect and the elements
of $\P(\pb^{(\la)},B)$ are called paths. The path realization of
$B(\la)$ is particularly well-suited to the study of branching
functions for the coset $\geh/\geh_{I\backslash\{0\}}$, where
$\geh_{I\backslash\{0\}}$ is the simple Lie subalgebra of $\geh$
whose Dynkin diagram is obtained from that of $\geh$ by removing
the vertex labeled $0$ in \cite{Kac}.

The purpose of this paper is to give a similar realization of a
finite tensor product of such highest weight modules.

Let $B_i$ ($i=1,2,\ldots,m$) be a perfect crystal of level $l_i$
of a quantum affine algebra $\Uq$ such that $l_1\ge l_2\ge
\cdots\ge l_m\ge l_{m+1}=0$, and $\la_i$ a dominant integral
weight of level $l_i-l_{i+1}$. We construct a set of paths
$\P(\pb^{(\la_1,\ldots,\la_m)},B_1\ot\cdots\ot B_m)$ with the
reference path $\pb^{(\la_1,\ldots,\la_m)}$ (section
\ref{subsec:paths}). For technical reasons we also require
assumptions (A1) and (A2) (section \ref{subsec:assumption}). These
assumptions are satisfied by many perfect crystals known so far,
such as $B_l$ for any nonexceptional affine Lie algebra $\geh$
given in \cite{KKM} (Example \ref{ex:KKM}) and $B^{k,l}$ for
$\geh=A^{(1)}_n$ given in \cite{KMN2} (Example \ref{ex:A}). With
this setup we have the following isomorphism of crystals (Theorem
\ref{th:main})
\begin{equation} \label{eq:main th}
\P(\pb^{(\la_1,\ldots,\la_m)},B_1\ot\cdots\ot B_m) \simeq
B(\la_1)\ot\cdots\ot B(\la_m).
\end{equation}
The main result of \cite{KMN1} is the case $m=1$.
A special case of $m=2$ for $\geh=A^{(1)}_n$ is treated in \cite{HKKOT}.

Our motivation for \eqref{eq:main th} comes from a conjecture
related to fermionic formulas given in \cite{HKOTT}. Fermionic
formulas are certain polynomials which originate from the
combinatorics of the Bethe Ansatz in solvable lattice models. They
are related to the representation theory of affine Lie algebras by
a conjecture in \cite{HKOTT} which states that they are equal to
one dimensional sums defined in terms of a conjectural family of
crystals of finite dimensional $\Uqp$-modules. If this series of
conjectures is shown to be true, we see from \eqref{eq:main th}
that fermionic formulas give truncated branching functions for the
coset $\geh/\geh_{I\backslash\{0\}}$ provided that the relevant
crystals are all perfect. For type $A^{(1)}_n$ the desired family
of crystals exists and all crystals therein are perfect
\cite{KMN2}, and the equality between the one dimensional sum and
the corresponding fermionic formula was established in \cite{KSS}.
Therefore the conjectures mentioned above are all settled in this
case.

To extend the result of \cite{KMN1} to the case of several tensor
factors, we employ some substantial results on affine crystals.
The first is the path realization of the crystal graph of the
lower triangular part $\Uq^-$ of $\Uq$ given in \cite{KKM} using
coherent families of perfect crystals. The second is the theory of
crystals with core \cite{KK}. Given a perfect crystal $B$ in a
coherent family and $B(\la)$, the crystal $B \otimes B(\la)$ is a
crystal with core, and there is a perfect crystal $B'$ and a
weight $\la'$ such that $B \otimes B(\la)\cong B(\la') \otimes
B'$. This result allows one to exchange a finite crystal past a
highest weight crystal. Given the above results, finding the
correct reference path and path space, more or less reduces to a
computation involving the basic isomorphisms of \cite{KMN1} and
\cite{KK} together with combinatorial $R$-matrices, which take the
form $B \otimes B' \cong B'\otimes B$ where $B$ and $B'$ are
perfect.

However the resulting isomorphism \eqref{eq:main th} so obtained,
is only known to preserve weights up to the null root $\delta$. As
in \cite{KMN1} one defines an energy function
$E:\P(\pb^{(\la_1,\dotsc,\la_m)},B_1\otimes\dotsm\otimes
B_m)\rightarrow \Z$ which determines the multiple of $\delta$ that
makes \eqref{eq:main th} weight-preserving. However unlike the
situation in \cite{KMN1}, our paths are inhomogeneous in the sense
that they have several different perfect crystals as tensor
factors. For inhomogeneous paths the evaluation of the energy
function requires the explicit computation of combinatorial
$R$-matrices, which are the identity map in the homogeneous case.
On minimal elements (defined in subsection \ref{ssec:finite}), it
appears that the combinatorial $R$-matrices can be computed
explicitly using automorphisms acting on perfect crystals. For
type $A^{(1)}_n$ this is proved in \cite{SS}. For general type we
must assume that such a result holds for the perfect crystals
being used. Several new and important observations regarding the
value of the energy function (in particular, its value on minimal
elements) are required to complete the proof of \eqref{eq:main
th}.

The plan of the paper is as follows. In section \ref{sec:crystal}
we prepare necessary notation and review crystals. The path space
$\P(\pb^{(\la_1,\dotsc,\la_m)},B_1\otimes\dotsm \otimes B_m)$ is
defined in section \ref{sec:paths}. In section \ref{sec:main} the
required assumptions are explained and the main theorem proved.
Section \ref{sec:A} is devoted to the demonstration that the
series of crystals $\{B^{k,l}\}_{l\ge1}$ of type $A^{(1)}_n$ given
in \cite{KMN2} forms a coherent family of perfect crystals for any
fixed $k$. This fact is necessary to prove the main theorem for
type $A^{(1)}_n$. We attach Appendix A for proofs of formulas of
$\et{0},\ft{0}$ used in section \ref{sec:A}.

\subsection*{Acknowledgments}
The authors thank Masaki Kashiwara for stimulating discussions.
M.O. thanks Goro Hatayama, Atsuo Kuniba and Taichiro Takagi for
collaboration in the early stage of this work.
M.O. is partially supported by Grant-in-Aid for Scientific Research from
the Ministry of Education, Culture, Sports, Science and Technology.
A.S. is partially supported by a Faculty Research Grant.
M.S. is partially supported by the grant NSF DMS-0100918.

\section{Crystals} \label{sec:crystal}

\subsection{Notation}
Let $\geh$ be an affine Lie algebra and $I$ the index set of its
Dynkin diagram. Note that $0$ is included in $I$. Let
$\alpha_i,h_i,\La_i$ ($i\in I$) be the simple roots, simple
coroots, and fundamental weights for $\geh$. Let
$\delta=\sum_{i\in I}a_i\alpha_i$ denote the standard null root
and $c=\sum_{i\in I}a_i^\vee h_i$ the canonical central element,
where $a_i,a_i^\vee$ are the positive integers given in
\cite{Kac}. Let $P=\bigoplus_{i\in I}\Z\La_i\oplus\Z\delta$ be the
weight lattice and $P^+=\sum_{i\in I}\Zn\La_i\oplus\Z\delta$ the
dominant weights.

Let $\Uq$ be the quantum affine algebra associated to $\geh$. For
the definition of $\Uq$ and its Hopf algebra structure, see e.g.
section 2.1 of \cite{KMN1}. For $J\subset I$ we denote by
$U_q(\geh_J)$ the subalgebra of $\Uq$ generated by $e_i,f_i,t_i$
($i\in J$). In particular, $U_q(\geh_{I\setminus\{0\}})$ is
identified with the quantized enveloping algebra for the simple
Lie algebra whose Dynkin diagram is obtained by deleting the
$0$-th vertex from that of $\geh$. We also consider the quantum
affine algebra without derivation $\Uqp$. The weight lattice of
$\geh'$ is called the classical weight lattice $\Pcl=P/\Z\delta$.
We canonically identify $\Pcl$ with $\bigoplus_{i\in
I}\Z\La_i\subset P$. For the precise treatment, see section 3.1 of
\cite{KMN1}. We further define the following subsets of $\Pcl$:
$\Pcl^0=\{\la\in \Pcl\mid \langle\la,c\rangle=0\}$,
$\Pcl^+=\{\la\in \Pcl\mid \langle\la,h_i\rangle\ge0\mbox{ for any
}i\}$, $\Pcll=\{\la\in \Pcl^+\mid \langle\la,c\rangle=l\}$. For
$\la,\mu\in\Pcl$, we write $\la\ge\mu$ to mean $\la-\mu\in\Pcl^+$.

\subsection{Crystals and crystal bases}
\label{ssec:crystal}
We summarize necessary facts in crystal theory. Our basic references are
\cite{K1}, \cite{KMN1} and \cite{AK}.

A crystal $B$ is a set $B=\sqcup_{\la\in P}B_\la$ ($\wt\, b=\la$
if $b\in B_\la$) with the maps
\begin{eqnarray*}
&&\et{i}: B_\la\longrightarrow B_{\la+\alpha_i}\sqcup\{0\},\quad
\ft{i}: B_\la\longrightarrow B_{\la-\alpha_i}\sqcup\{0\},\\
&&\veps_i : B\longrightarrow \Z\sqcup\{-\infty\},\quad
\vphi_i : B\longrightarrow \Z\sqcup\{-\infty\}
\end{eqnarray*}
for all $i\in I$ such that
\begin{eqnarray}
&&\mbox{for $b\in B_\la$, $\vphi_i(b)=\langle h_i,\la\rangle+\veps_i(b)$},
\label{eq:i-wt}\\
&&\mbox{for $b\in B$, we have}\\
&&\hspace{1.5cm}
\begin{array}{rcl}
\veps_i(b)&=&\veps_i(\et{i}b)+1\mbox{ if }\et{i}b\neq0,\\
&=&\veps_i(\ft{i}b)-1\mbox{ if }\ft{i}b\neq0,\\
\vphi_i(b)&=&\vphi_i(\et{i}b)-1\mbox{ if }\et{i}b\neq0,\\
&=&\vphi_i(\ft{i}b)+1\mbox{ if }\ft{i}b\neq0,
\end{array} \nonumber\\
&&\mbox{for $b,b'\in B$, $\et{i}b'=b$ if and only if $b'=\ft{i}b$},\\
&&\mbox{for $b\in B$, $\veps_i(b)=\vphi_i(b)=-\infty$ implies
$\et{i}b=\ft{i}b=0$}.
\end{eqnarray}
A crystal $B$ can be regarded as a colored oriented graph by defining
\[
b\stackrel{i}{\longrightarrow}b'\quad\Longleftrightarrow\quad \ft{i}b=b'.
\]
If we want to emphasize $I$, $B$ is called an $I$-crystal.

Important examples of crystals are given by the crystal bases of
integrable $\Uq$ (or $\Uqp$)-modules. They satisfy
\begin{itemize}
\item[] for any $b$ and $i$, there exists $n>0$ such that
$\et{i}^nb=\ft{i}^nb=0$.
\end{itemize}
In such cases the maps $\veps_i,\vphi_i$ are given by
\[
\veps_i(b)=\max\{n\in\Zn\mid\et{i}^nb\neq0\},\quad
\vphi_i(b)=\max\{n\in\Zn\mid\ft{i}^nb\neq0\}.
\]
We also set
\[
\veps(b)=\sum_{i\in I}\veps_i(b)\La_i,\quad
\vphi(b)=\sum_{i\in I}\vphi_i(b)\La_i.
\]
If $B$ is $\Pcl$-weighted, i.e., $\wt b\in\Pcl$ for any $b\in B$,
(\ref{eq:i-wt}) is equivalent to $\vphi(b)-\veps(b)=\wt b$.

For two crystals $B_1$ and $B_2$ a {\em morphism} of crystals from $B_1$
to $B_2$ is a map $\psi : B_1\sqcup\{0\}\longrightarrow B_2\sqcup\{0\}$
such that
\begin{eqnarray}
&&\psi(0)=0,\\
&&\psi(\et{i}b)=\et{i}\psi(b)\mbox{ for }b,\et{i}b\in B_1 \mbox{ and }
\psi(\ft{i}b)=\ft{i}\psi(b)\mbox{ for }b,\ft{i}b\in B_1, \label{eq:morph ef}\\
&&\mbox{for }b\in B_1,\veps_i(b)=\veps_i(\psi(b)),\vphi_i(b)=\vphi_i(\psi(b))
\mbox{ if }\psi(b)\in B_2,\\
&&\mbox{for }b\in B_1,\wt b=\wt \psi(b)\mbox{ if }\psi(b)\in B_2.
\end{eqnarray}
A morphism of crystals $\psi : B_1\longrightarrow B_2$ is called an
{\em embedding} if $\psi$ is injective.

For two crystals $B_1$ and $B_2$, the tensor product $B_1\ot B_2$
is defined.
\[
B_1\ot B_2=\{b_1\ot b_2\mid b_1\in B_1,b_2\in B_2\}.
\]
The actions of $\et{i}$ and $\ft{i}$ are defined by
\begin{eqnarray}
  \label{eq:ot-e} \et{i}(b_1\ot b_2)&=&\begin{cases}
  \et{i}b_1\ot b_2 & \text{if $\vphi_i(b_1)\ge\veps_i(b_2)$} \\
  b_1\ot \et{i}b_2 & \text{if $\vphi_i(b_1) < \veps_i(b_2)$,}
  \end{cases} \\
  \label{eq:ot-f}
  \ft{i}(b_1\ot b_2)&=& \begin{cases}
  \ft{i}b_1\ot b_2 & \text{if $\vphi_i(b_1) > \veps_i(b_2)$} \\
  b_1\ot \ft{i}b_2 & \text{if $\vphi_i(b_1)\le\veps_i(b_2)$.}
  \end{cases}
\end{eqnarray}
Here $0\ot b$ and $b\ot0$ are understood to be $0$.
$\veps_i,\vphi_i$ and $\wt$ are given by
\begin{eqnarray}
\veps_i(b_1\ot b_2)&=&
\max(\veps_i(b_1),\veps_i(b_1)+\veps_i(b_2)-\vphi_i(b_1)),\label{eq:ot-eps}\\
\vphi_i(b_1\ot b_2)&=&
\max(\vphi_i(b_2),\vphi_i(b_1)+\vphi_i(b_2)-\veps_i(b_2)),\label{eq:ot-phi}\\
\wt(b_1\ot b_2)&=&\wt b_1+\wt b_2.
\end{eqnarray}
With this tensor product operation, $I$-crystals form a tensor
category.

The following crystal $T_\la$ will be used later.

\begin{example} \label{ex:T}
For $\la\in\Pcl$ consider the set $T_\la=\{t_\la\}$ with one element.
Set $\veps_i(t_\la)=\vphi_i(t_\la)=-\infty$ for $i\in I$ and $\wt t_\la=\la$.
Given a crystal $B$ one can take the tensor product
\[
T_\la\ot B\ot T_\mu=\{t_\la\ot b\ot t_\mu\mid b\in B\}.
\]
Then we have
\begin{alignat*}{3}
\wt(t_\la\ot b\ot t_\mu)&=\la+\mu+\wt b, & & \\
\et{i}(t_\la\ot b\ot t_\mu)&=t_\la\ot\et{i}b\ot t_\mu, &\quad
\veps_i(t_\la\ot b\ot t_\mu)&=\veps_i(b)-\langle
h_i,\la\rangle \\
\ft{i}(t_\la\ot b\ot t_\mu)&=t_\la\ot\ft{i}b\ot t_\mu, & \quad
\vphi_i(t_\la\ot b\ot t_\mu)&=\vphi_i(b)+\langle h_i,\mu\rangle.
\end{alignat*}
\end{example}

\begin{definition}[\cite{AK}]
We say a $P$ (or $\Pcl$)-weighted crystal is regular, if for any $i,j\in I$
($i\ne j$), $B$ regarded as $\{i,j\}$-crystal is a disjoint union of
crystals of integrable highest weight modules over $U_q(\geh_{\{i,j\}})$.
\end{definition}

Let $V(\la)$ be the integrable highest weight $\Uq$-module
with highest weight $\la\in P^+$ and highest weight vector $u_\la$. It is
shown in \cite{K1} that $V(\la)$ has a crystal base $(L(\la),B(\la))$.
We regard $u_\la$ as an element of $B(\la)$ as well. $B(\la)$ is a regular
$P$-weighted crystal. A finite-dimensional integrable $\Uqp$-module $V$
does not necessarily have a crystal base. If $V$ has a crystal base $(L,B)$,
then $B$ is a regular $\Pcl^0$-weighted crystal with finitely many
elements.

Let $W$ be the affine Weyl group associated to $\geh$, and $s_i$ be the
simple reflection corresponding to $\alpha_i$. $W$ acts on any regular
crystal $B$ \cite{K2}. The action is given by
\[
S_{s_i}b=\left\{
\begin{array}{ll}
\ft{i}^{\langle h_i,\wts b\rangle}b\quad
&\mbox{ if }\langle h_i,\wt b\rangle\ge0\\
\et{i}^{-\langle h_i,\wts b\rangle}b\quad
&\mbox{ if }\langle h_i,\wt b\rangle\le0.
\end{array}
\right.
\]
An element $b$ of $B$ is called $i$-{\em extremal} if $\et{i}b=0$ or
$\ft{i}b=0$. $b$ is called {\em extremal} if $S_wb$ is $i$-extremal
for any $w\in W$ and $i\in I$.

\begin{definition}[\cite{AK} Definition 1.7] \label{def:simple}
Let $B$ be a regular $\Pcl^0$-weighted crystal with finitely many
elements.
We say $B$ is simple if it satisfies
\begin{itemize}
\item[(1)] There exists $\la\in\Pcl^0$ such that the weights of $B$ are
in the convex hull of $W\la$.
\item[(2)] $\sharp B_\la=1$.
\item[(3)] The weight of any extremal element is in $W\la$.
\end{itemize}
\end{definition}

Acting by the Weyl group of the canonical simple Lie subalgebra,
one may choose the above weight $\la$ such that
\begin{equation} \label{eq:Jdom}
\langle h_i\,,\,\la\rangle\ge0 \qquad\text{ for $i\in I\backslash\{0\}$. }
\end{equation}
For $B$ simple, let $u(B)\in B$ be the unique element of weight $\la$
such that $\la$ satisfies \eqref{eq:Jdom}.

\begin{prop}[\cite{AK} Lemma 1.9 \& 1.10] \label{prop:simple}
Simple crystals have the following properties.
\begin{itemize}
\item[(1)] A simple crystal is connected.
\item[(2)] The tensor product of simple crystals is also simple.
\end{itemize}
\end{prop}

\subsection{Category $\Cfin$}
\label{ssec:finite} We recall the category $\Cfin$ defined in
\cite{HKKOT}. Let $B$ be a regular $\Pcl^0$-weighted crystal with
finitely many elements. For $B$ we introduce the {\em level} of
$B$ by
\[
\lev B=\min\{\langle c,\veps(b)\rangle\mid b\in B\}\in\Zn.
\]
Note that $\langle c,\veps(b)\rangle=\langle c,\vphi(b)\rangle$ for any
$b\in B$. We also set $\Bmin=\{b\in B\mid\langle c,\veps(b)\rangle=\lev B\}$
and call an element of $\Bmin$ {\em minimal}.

\begin{definition}[\cite{HKKOT} Definition 2.5] \label{def:C-fin}
We denote by $\Cfin(\geh)$ (or simply $\Cfin$) the category of crystals
$B$ satisfying the following conditions:
\begin{itemize}
\item[(1)] $B$ is a crystal base of a finite-dimensional $\Uqp$-module.
\item[(2)] $B$ is simple.
\item[(3)] For any $\la\in\Pcl^+$ such that $\langle c,\la\rangle\ge
           \lev B$, there exists $b\in B$ satisfying $\veps(b)\le\la$.
           It is also true for $\vphi$.
\end{itemize}
\end{definition}

We call an element of $\Cfin(\geh)$ a {\em finite crystal}.

\begin{remark}\mbox{}
\begin{itemize}
\item[(i)] Condition (1) implies that $B$ is a regular $\Pcl^0$-weighted
           crystal with finitely many elements.
\item[(ii)] Set $l=\lev B$. Condition (3) implies that the maps $\veps$
           and $\vphi$ from $\Bmin$ to $\Pcll$ are surjective
           (cf. (4.6.5) in \cite{KMN1}).
\item[(iii)] Practically, one has to check condition (3) only for
           $\la\in\Pcl^+$ such that there is no $i\in I$ satisfying
           $\la-\La_i\ge0$ and $\langle c,\la-\La_i\rangle\ge\lev B$.
           In particular, if $a^\vee_i=1$ for any $i\in I$,
           that is, when $\geh=A_n^{(1)},C_n^{(1)}$,
           the surjectivity of $\veps$ and $\vphi$ assures (3).
\item[(iv)] The authors do not know a crystal satisfying (1) and (2), but
           not (3).
\end{itemize}
\end{remark}

Let $B_1$ and $B_2$ be two finite crystals. Definition
\ref{def:C-fin} (1) and the existence of the universal $R$-matrix assures
that there is a natural isomorphism of crystals
\begin{equation} \label{eq:iso}
B_1\ot B_2\simeq B_2\ot B_1.
\end{equation}

The following lemma is immediate.

\begin{lemma} \label{lem:tensor}
Let $B_1,B_2$ be finite crystals.
\begin{itemize}
\item[(1)] $\lev(B_1\ot B_2)=\max(\lev B_1,\lev B_2)$.
\item[(2)] If $\lev B_1\ge\lev B_2$, then $(B_1\ot B_2)_{\min}=
           \{b_1\ot b_2\mid b_1\in(B_1)_{\min},$
           $\vphi(b_1)\ge\veps(b_2)\}$.
\item[(3)] If $\lev B_1\le\lev B_2$, then $(B_1\ot B_2)_{\min}=
           \{b_1\ot b_2\mid b_2\in(B_2)_{\min},$
           $\vphi(b_1)\le\veps(b_2)\}$.
\end{itemize}
\end{lemma}

$\Cfin(\geh)$ is a tensor category, as it is closed under the
tensor product operation.

\subsection{Perfect crystals} \label{subsec:perfect}
An object $B$ of $\Cfin(\geh)$ is called {\em perfect} if the maps
$\veps$ and $\vphi$ are bijections from $\Bmin$ onto $\Pcll$,
where $l=\lev B$. If $B$ is perfect of level $l$, there is a
bijection $\sigma:\Pcll\rightarrow\Pcll$ defined by
$\sigma=\veps\circ\vphi^{-1}$. It is called the associated
automorphism of $B$. This nomenclature is explained below.

The next proposition is immediate from Lemma \ref{lem:tensor}.

\begin{prop}
Let $B$ be a perfect crystal of level $l$ with associated
automorphism $\sigma$. Then $B^{\ot m}$ is perfect of the same
level with associated automorphism $\sigma^m$.
\end{prop}

For a finite crystal $B$ write $B^{\le\la}=\{b\in B\mid \veps(b)\le\la\}$.

\begin{theorem}[\cite{KMN1}] \label{th:sum-of-B}
Assume $\mbox{\sl rank }\geh>2$.
Let $B$ be a perfect crystal of level $l$. For any $\la\in(\Pcl^+)_k$ with
$k\ge l$, there exists an isomorphism
\[
B(\la)\ot B\simeq\bigoplus_{b\in B^{\le\la}}B(\la+\wt b)
\]
of $\Pcl$-weighted crystals. In particular, if $k=l$, we have
\[
B(\la)\ot B\simeq B(\sigma^{-1}\la),
\]
where $\sigma$ is the associated automorphism of $B$.
\end{theorem}

We now recall the notion of a coherent family of perfect crystals
\cite{KKM}. Let $\{B^l\}_{l\ge1}$ be a family of perfect crystals
$B^l$ of level $l$. Take the index set $J=\{(l,b)\mid
l\in\Z_{>0},b\in B^l_{\min}\}$. Let $T_\la$ be the crystal in
Example \ref{ex:T}. A crystal $B^\infty$ with distinguished
element $b_\infty$ is called a {\em limit} of $\{B^l\}_{l\ge1}$ if
it satisfies the following conditions:
\begin{eqnarray}
&&\wt b_\infty=0,\veps(b_\infty)=\vphi(b_\infty)=0,\\
&&\mbox{for any $(l,b)\in J$, there exists an embedding of crystals}\\
&&\hspace{1.5cm}f_{(l,b)} : T_{\veps(b)}\ot B^l\ot T_{-\vphi(b)}\longrightarrow
B^\infty \nonumber\\
&&\mbox{sending $t_{\veps(b)}\ot b\ot t_{-\vphi(b)}$ to $b_\infty$,} \nonumber\\
&&B^\infty=\bigcup_{(l,b)\in J}\mbox{Im }f_{(l,b)}.
\end{eqnarray}

If a limit exists for the family $\{B^l\}$, we say $\{B^l\}$ is a
{\em coherent family} of perfect crystals. Set
$B^\infty_{\min}=\{b\in B^\infty\mid\langle
c,\veps(b)\rangle=0\}$. Then both $\veps$ and $\vphi$ map
$B^\infty_{\min}$ bijectively to $\Pcl^0$. It follows from
\cite[Lemma 4.6]{KKM} that the bijection given by
$\sigma=\veps\circ\vphi^{-1}$, is a linear automorphism of
$\Pcl^0$.

\begin{conjecture} \label{conj:auto} \cite{KK}
\begin{enumerate}
\item  For any coherent family of perfect crystals,
\begin{equation}
\begin{split}
&\text{$\sigma$ extends to a linear automorphism of $\Pcl$}\\
&\text{such that $\sigma\vphi(b)=\veps(b)$ for any $b\in
B^l_{\min}$.} \label{eq:sigma-cond}
\end{split}
\end{equation}
\item $\sigma$ is induced by a Dynkin diagram automorphism.
That is, there is an automorphism $\tau:I\rightarrow I$ of the
Dynkin diagram of $\geh$ such that $\sigma(\La_i)=\La_{\tau(i)}$
for all $i\in I$.
\end{enumerate}
\end{conjecture}

\begin{lemma} \label{lem:auto} Let $\{B^l\}$ be a coherent
family of perfect crystals. Suppose there is an $i_0\in I$ such
that $\langle c,\La_{i_0}\rangle=1$ and $\wt
\vphi^{-1}_{B^l}(l\La_{i_0})=l\cdot\wt \vphi^{-1}_{B^1}(\La_{i_0})$ for
all $l\in\mathbb{Z}_{>0}$. Then \eqref{eq:sigma-cond} holds.
\end{lemma}
\begin{proof} Let $v_l=\vphi^{-1}_{B^l}(l\La_{i_0})$ for $l\in\mathbb{Z}_{>0}$.
The hypotheses imply that
$\sigma_{B^l}(l\La_{i_0})=l\sigma_{B^1}(\La_{i_0})$. It follows
that there is a unique well-defined linear automorphism $\sigma$
of $\Pcl$ such that $\sigma(\La_{i_0})=\sigma_{B^1}(\La_{i_0})$
and $\sigma|_{\Pcl^0}=\veps\circ \vphi^{-1}$. Let $\la\in \Pcll$
and $b\in B^l_{\min}$ such that $\vphi_{B^l}(b)=\la$. Then
$\la-l\La_{i_0}\in \Pcl^0$. Write
$b'=f_{(l,v_l)}(t_{\veps(v_l)}\otimes b\otimes t_{-l\La_{i_0}})$.
By the hypotheses we have
\begin{equation*}
\begin{split}
  0&=\sigma \vphi(b')-\veps(b') \\
  &=
  \sigma(\vphi_{B^l}(b)-l\La_{i_0})-(\veps_{B^l}(b)-\sigma_{B^l}(l\La_{i_0}))\\
  &=
  \sigma(\la)-l\sigma(\La_{i_0})-\sigma_{B^l}(\la)+l\sigma_{B^1}(\La_{i_0})\\
  &=\sigma(\la)-\sigma_{B^l}(\la).
\end{split}
\end{equation*}
This is precisely \eqref{eq:sigma-cond}.
\end{proof}

\begin{lemma} \label{lem:Dynkin} Suppose \eqref{eq:sigma-cond}
holds. Then there is a permutation $\tau:I\rightarrow I$ such that
$\sigma(\La_i)=\La_{\tau(i)}$ for all $i\in I$.
\end{lemma}
\begin{proof} By assumption $\sigma=\sigma_{B^l}$ on $\Pcll$. Thus $\sigma$
permutes $\Pcll$ for $l>0$. By linearity it follows that $\sigma$
permutes the dominant weights of level $l$ that cannot be
expressed as the sum of two nonzero dominant weights of smaller
level. But these are precisely the fundamental weights of level
$l$.
\end{proof}

Note that this does not prove that $\tau$ must preserve the Dynkin
diagram.

\subsection{Crystals with core}
In \cite{KK} Kang and Kashiwara developed the theory of {\em crystals with core}.
Let $B$ be a regular crystal. For $b\in B$ define
\[
{\mathcal E}(b)
=\{\et{i_l}\cdots\et{i_1}b\mid l\ge0,i_1,\ldots,i_l\in I\}\setminus\{0\}.
\]

\begin{definition}[\cite{KK} Definition 3.4]
We say a regular crystal $B$ has a core if ${\mathcal E}(b)$ is a finite set
for any $b\in B$. In this case, we define the core $C(B)$ of $B$ by
\[
C(B)=\{b\in B\mid {\mathcal E}(b')={\mathcal E}(b)
\mbox{ for every }b'\in {\mathcal E}(b)\},
\]
and $B$ is called a crystal with core.
\end{definition}

The following results are proven in \cite{KK}.

\begin{theorem}[\cite{KK} Theorem 5.4] \label{th:core}
Suppose $\mbox{\sl rank }\geh>2$. Let $\{B^l\}_{l\ge1}$ be a coherent family of
perfect crystals with the associated automorphism $\sigma$ satisfying
(\ref{eq:sigma-cond}). Take positive integers $k,l$ such that $k<l$ and
a weight $\la\in (\Pcl^+)_k$. Then we have an isomorphism of $\Pcl$-weighted
crystals
\begin{equation} \label{eq:iso-with-core}
B(\la)\ot B^l\simeq B^{l-k}\ot B(\sigma^{-1}\la).
\end{equation}
\end{theorem}

\begin{prop}[\cite{KK}] \label{prop:core}
On the core of both sides of (\ref{eq:iso-with-core}) we have
\begin{itemize}
\item[(1)] $C(B^{l-k}\ot B(\sigma^{-1}\la))=B^{l-k}\ot u_{\sigma^{-1}\la}$,
\item[(2)] $C(B(\la)\ot B^l)=u_\la\ot (B^l)^{(\la)}$ with a suitable subset
           $(B^l)^{(\la)}$ of $B^l$. Moreover, $(B^l)^{(\la)}\cap B^l_{\min}
            =\{b\in B^l\mid\veps(b)-\la\in(\Pcl^+)_{l-k}\}$.
\item[(3)] For any distinct elements $b\in B(\la)\ot B^l$ and $b'\in
           C(B(\la)\ot B^l)$, there exists a sequence $i_1,\ldots,i_l\in I$
           ($l\ge1$) such that $b'\in\et{i_l}\cdots\et{i_1}b$.
\end{itemize}
\end{prop}

\section{Paths} \label{sec:paths}

In this section we review the set of paths $\PB$ defined in \cite{HKKOT}
and prepare necessary facts.

\subsection{Energy function}
Let $B_1$ and $B_2$ be finite crystals.
Suppose $b_1\ot b_2\in B_1\ot B_2$ is mapped to $\bt_2\ot\bt_1
\in B_2\ot B_1$ under the isomorphism (\ref{eq:iso}). A $\Z$-valued function
$H$ on $B_1\ot B_2$ is called an {\em energy function} if for any $i$
and $b_1\ot b_2\in B_1\ot B_2$ such that $\et{i}(b_1\ot b_2)\neq0$,
 it satisfies
\begin{equation} \label{eq:e-func}
H(\et{i}(b_1\ot b_2))=
\begin{cases}
H(b_1\ot b_2)+1
&\mbox{ if }i=0,\vphi_0(b_1)\geq\veps_0(b_2),\vphi_0(\bt_2)\geq\veps_0(\bt_1),\\
H(b_1\ot b_2)-1
&\mbox{ if }i=0,\vphi_0(b_1)<\veps_0(b_2),\vphi_0(\bt_2)<\veps_0(\bt_1),\\
H(b_1\ot b_2)&\mbox{ otherwise}.
\end{cases}
\end{equation}
When we want to emphasize $B_1\ot B_2$, we write $H_{B_1B_2}$ for $H$.
The existence of such a function can be shown in a similar manner to
section 4 of \cite{KMN1} based on the existence of the {\em combinatorial
$R$-matrix}. The energy function is unique up to an additive constant, since
$B_1\ot B_2$ is connected. By definition, $H_{B_1B_2}(b_1\ot b_2)
=H_{B_2B_1}(\bt_2\ot\bt_1)$.

If the tensor product $B_1\ot B_2$ is homogeneous, i.e., $B_1=B_2$,
we have $\bt_2=b_1,\bt_1=b_2$. Thus (\ref{eq:e-func}) reduces to
\begin{equation} \label{eq:e-func*}
H(\et{i}(b_1\ot b_2))=
\begin{cases}
H(b_1\ot b_2)+1
&\mbox{ if }i=0,\vphi_0(b_1)\geq\veps_0(b_2),\\
H(b_1\ot b_2)-1
&\mbox{ if }i=0,\vphi_0(b_1)<\veps_0(b_2),\\
H(b_1\ot b_2)&\mbox{ if }i\ne0.
\end{cases}
\end{equation}

Let $M=[i_1,\ldots,i_l]$ be a multiset from $I$. A partition of $M$ is a decomposition
$M=[j_1,\ldots,j_m]\sqcup [j_1^\ast,\ldots,j_{m^\ast}^\ast]$ such that
$[i_1,\ldots,i_l]=[j_1,\ldots,j_m,j_1^\ast,\ldots,j_{m^\ast}^\ast]$.
For $r\in I$ and a multiset $M=[i_1\ldots,i_l]$ from $I$ let
$\sharp r[i_1,\ldots,i_l]$ denote the number of letters $r$ in $M$.

The next lemma is an easy consequence of the definition of $H$.

\begin{lemma} \label{lem:diff-H}
Suppose $b_1\ot b_2$ is mapped to $\bt_2\ot \bt_1$ under the isomorphism
$B_1\ot B_2\simarrow B_2\ot B_1$
and $b'_1\ot b'_2=\et{i_l}\cdots\et{i_1}(b_1\ot b_2)$.
If we have
\begin{eqnarray*}
\et{i_l}\cdots\et{i_1}(b_1\ot b_2)&=&
\et{j_{m^\ast}^\ast}\cdots\et{j_1^\ast}b_1\ot\et{j_m}\cdots\et{j_1}b_2\\
\mbox{and}\quad
\et{i_l}\cdots\et{i_1}(\bt_2\ot\bt_1)&=&
\et{j'_{m'}}\cdots\et{j'_1}\bt_2\ot
\et{j^{\prime\ast}_{m^{\prime\ast}}}\cdots\et{j^{\prime\ast}_1}\bt_1,
\end{eqnarray*}
then
\[
H(b'_1\ot b'_2)-H(b_1\ot b_2)
=\sharp0[j'_1,\ldots,j'_{m'}]-\sharp0[j_1,\ldots,j_m].
\]
\end{lemma}

The following proposition reduces the energy function of tensor products
to that of their components. See section 2.13 of \cite{OSS}.

\begin{prop} \label{prop:HofTensor}
Consider the tensor product $B=B_1\ot B_2\ot\cdots\ot B_m$.
\begin{itemize}
\item[(1)] Set $B^\star=B_1\ot\cdots\ot B_{m-1}$. For $b_1\ot\cdots\ot b_m\in B$,
define $b_m^{(j)}$ ($1\le j\le m$) by
\begin{equation} \label{eq:b^(j)}
\begin{array}{ccc}
B_j\ot\cdots\ot B_{m-1}\ot B_m&\simarrow&B_m\ot B_j\ot\cdots\ot B_{m-1}\\
b_j\ot\cdots\ot b_{m-1}\ot b_m&\mapsto&b_m^{(j)}\ot\bt_j\ot\cdots\ot\bt_{m-1}.
\end{array}
\end{equation}
We understand $b_m^{(m)}=b_m$.
Then we have
\[
H_{B^\star B_m}((b_1\ot\cdots\ot b_{m-1})\ot b_m)=\sum_{1\le j\le m-1}H_{B_jB_m}(b_j\ot b^{(j+1)}_m).
\]

\item[(2)] Set $B^\dagger=B_2\ot\cdots\ot B_m$. For $b_1\ot\cdots\ot b_m\in B$,
define $b_1^{\langle j\rangle}$ ($1\le j\le m$) by
\begin{equation} \label{eq:b^<j>}
\begin{array}{ccc}
B_1\ot B_2\ot\cdots\ot B_j&\simarrow&B_2\ot\cdots\ot B_j\ot B_1\\
b_1\ot b_2\ot\cdots\ot b_j&\mapsto&\hat{b}_2\ot\cdots\ot\hat{b}_j\ot b_1^{\langle j\rangle}.
\end{array}
\end{equation}
We understand $b_1^{\langle 1\rangle}=b_1$.
Then we have
\[
H_{B_1B^\dagger}(b_1\ot(b_2\ot\cdots\ot b_m))=\sum_{2\le j\le m}H_{B_1B_j}(b_1^{\langle j-1\rangle}\ot b_j).
\]

\end{itemize}
\end{prop}

\subsection{Energy $D_B$} We wish to define an energy function
$D_B:B\rightarrow \Z$ for tensor products of perfect crystals
of the form $B^{r,s}$.

Let $B=B^{r,s}$ be perfect. Then there exists
a unique element $b^\na\in B$ such that $\varphi(b^\na)=\lev(B)\La_0$.
Define $D_B:B\rightarrow\Z$ by
\begin{equation*}
  D_B(b) = H_{BB}(b^\na\otimes b) - H_{BB}(b^\na \otimes u(B))
\end{equation*}
where $u(B)$ is defined in subsection \ref{ssec:crystal}.
From now on we assume that the energy function
is normalized by
\begin{equation} \label{eq:Hnorm}
H_{B_1B_2}(u(B_1)\otimes u(B_2))=0.
\end{equation}

Now suppose by induction that $B_i$ is a tensor product of
perfect crystals of the form $B^{r,s}$ such that $D_{B_i}$ has been
defined for $1\le i\le m$. For $B=B_1\otimes\dotsm\otimes B_m$
and $b_i\in B_i$, define
\begin{equation} \label{eq:func}
  D_B(b_1\ot\dotsm\ot b_m) = \sum_{1\le i\le m} D_{B_i}(b_i^{(1)}) +
  \sum_{1\le i<j\le m} H_{B_iB_j}(b_i\otimes b_j^{(i+1)})
\end{equation}
where $b_j^{(i)}$ is defined by \eqref{eq:b^(j)}.

\begin{prop}[\cite{OSS} Proposition 2.13] \label{pp:Ddef}  The above
function $D_B$ is well-defined, that is, it depends
only on $B$ and not on the grouping of tensor factors within $B$.
\end{prop}

\begin{example} Let $B_i$ have the form $B^{r,s}$ for $1\le i\le 4$
and let $B=B_1\otimes B_2\otimes B_3\otimes B_4$. One may define
$D_B$ using three groups of tensor factors
by the grouping $(B_1) \otimes (B_2\otimes B_3) \otimes (B_4)$
or using two tensor factors by the grouping
$(B_1\otimes B_2) \otimes (B_3\otimes B_4)$. Proposition \ref{pp:Ddef}
asserts that the resulting functions $D_B$ are equal.
\end{example}

\begin{prop}[\cite{OSS} Proposition 2.15] \label{prop:D intrinsic}
Let $(i_1,i_2,\ldots,i_m)$ be a permutation of the set $\{1,2,\ldots,m\}$, and
set $\tilde{B}=B_{i_1}\ot B_{i_2}\ot\cdots\ot B_{i_m}$. Suppose that $b\in B$ is
mapped to $\bt\in\tilde{B}$ under the isomorphism $B\simarrow\tilde{B}$. Then
we have $D_B(b)=D_{\tilde{B}}(\bt)$.
\end{prop}

If the tensor product $B_1\ot\cdots\ot B_m$ is homogeneous, i.e.,
$B_1=\cdots=B_m$, we have $b^{(i+1)}_j=b_{i+1}$
for any $1\le i<j\le m$ and $b_i^{(1)}=b_1$ for $1\le i\le m$.
Then \eqref{eq:func} reduces to
\begin{equation} \label{eq:hom Denergy}
D_B(b)=m D_{B_1}(b_1) + \sum_{j=1}^{m-1}(m-j)H_{B_1B_1}(b_j\ot b_{j+1}).
\end{equation}
For later use we prepare a lemma, which is an immediate consequence of
\eqref{eq:e-func*}.

\begin{lemma} \label{lem:D diff}
Suppose $B=B_1^{\ot m}$. For $b=b_1\ot\cdots\ot b_m\in B$ and $i\in I$
such that $\et{i}b\neq0$ we have
\[
D_B(\et{i}b)=D_B(b)-\delta_{i0},
\]
unless $i=0$ and $\et{0}b=\et{0}b_1\ot b_2\ot\cdots\ot b_m$.
\end{lemma}

\subsection{Set of paths $\PB$} \label{subsec:paths}
An element $\pb=\cdots\ot\bb_j\ot\cdots\ot\bb_2\ot\bb_1$ of the semi-infinite
tensor product of $B$ is called a {\em reference path}
if it satisfies $\bb_j\in \Bmin$ and
$\vphi(\bb_{j+1})=\veps(\bb_j)$ for any $j\ge1$.
Fix a reference path $\pb=\cdots\ot\bb_j\ot\cdots\ot\bb_2\ot\bb_1$.
The {\em set of paths} $\PB$ is defined by
\[
\PB=\{p=\cdots\ot b_j\ot\cdots\ot b_2\ot b_1\mid b_j\in B,b_k=\bb_k\mbox{ for }
k\gg1\}.
\]
An element of $\PB$ is called a {\em path}. For convenience we denote
$b_k$ by $p(k)$ and $\cdots\ot b_{k+2}\ot b_{k+1}$ by $p[k]$ for
$p=\cdots\ot b_j\ot\cdots\ot b_2\ot b_1$.
For a path $p\in\PB$, set
\begin{eqnarray}
E(p)&=&\sum_{j=1}^\infty j(H(p(j+1)\ot p(j))-H(\pb(j+1)\ot\pb(j))),
\label{eq:def E}\\
\wtb p&=&\vphi(\pb(1))+\sum_{j=1}^\infty (\wt p(j)-\wt \pb(j))
-(E(p)/a_0)\delta, \label{eq:bold wt}
\end{eqnarray}
where $a_0$ is the $0$-th Kac label.
$E(p)$ and $\wtb p$ are called the {\em energy} and {\em weight} of $p$.
We distinguish $\wtb p\in P$ from
$\wt p=\vphi(\pb(1))+\sum_{j=1}^\infty (\wt p(j)-\wt \pb(j))\in\Pcl$.
Compare $E(p)$ with \eqref{eq:hom Denergy}. It is a normalized energy for
the semi-infinite tensor product of $B$.

\begin{remark}\mbox{}
\begin{itemize}
\item[(i)] If $B$ is perfect, the set of reference paths is in bijection
           with $\Pcll$, where $l=\lev B$. For $\la\in\Pcll$ take a unique
           $\bb_1\in\Bmin$ such that $\vphi(\bb_1)=\la$. The condition
           $\vphi(\bb_{j+1})=\veps(\bb_j)$ fixes
           $\pb=\cdots\ot\bb_j\ot\cdots\ot\bb_1$ uniquely.
\item[(ii)] In \cite{KMN1} $\pb$ is called a ground state path,
           since $E(p)\ge E(\pb)$ for any $p\in\PB$. But if $B$ is not
           perfect, it is no longer true in general.
\end{itemize}
\end{remark}

The following theorem is important.

\begin{theorem}[\cite{HKKOT} Theorem 3.7] \label{th:C-h}
Assume $\mbox{\sl rank }\geh>2$. Then $\PB$ is isomorphic to
a direct sum of crystals $B(\la)$ ($\la\in P^+$) of integrable
highest weight $\Uq$-modules.
\end{theorem}

For a set of paths $\PB$ let $\PB_0$ denote the set of its
{\em highest weight} elements, i.e., elements satisfying $\et{i}p=0$
for any $i$.
The following proposition describes the set of highest weight elements in
$\PB$.

\begin{prop}[\cite{HKKOT} Proposition 3.9 \& Corollary 3.10] \label{prop:P_0}
\[
\PB_0=\{p\in\PB\mid p(j)\in\Bmin,
\vphi(p(j+1))=\veps(p(j))\mbox{ for }\forall j\}.
\]
Moreover, if $p\in\PB_0$, then $\wt p[j]=\vphi(p(j+1))$.
\end{prop}

{}From Theorem \ref{th:C-h} and Proposition \ref{prop:P_0} we obtain

\begin{theorem}[\cite{KMN1} Proposition 4.6.4] \label{th:KMN}
Let $B$ be a perfect crystal of level $l$ with the associated automorphism $\sigma$.
For $\la\in(\Pcl^+)_l$ define a reference path $\pb^{(\la)}$ by $\pb^{(\la)}(j)=\veps_B^{-1}(\sigma^j\la)$.
Then we have an isomorphism of $P$-weighted crystals
\[
\P(\pb^{(\la)},B)\simeq B(\la).
\]
\end{theorem}

\section{Tensor product theorem} \label{sec:main}

In this section we assume $\mbox{\sl rank }\geh>2$ because of frequent use
of Theorems \ref{th:sum-of-B}, \ref{th:core} and \ref{th:C-h}.

\subsection{Assumptions} \label{subsec:assumption}
In this subsection we explain two properties (A1) and (A2) for
a set of perfect crystals that are assumed to prove our main theorem.
Consider a set of perfect crystals
\[
{\mathcal B}=\{B_1,B_2,\ldots,B_m\}.
\]
We assume $\lev B_1\ge\lev B_2\ge\cdots\ge\lev B_m$.

\begin{itemize}
\item[(A1)]
Each $B_i$ in ${\mathcal B}$ belongs to a certain coherent family satisfying
the condition (\ref{eq:sigma-cond}).
\end{itemize}

\noindent
It should be emphasized that the coherent family containing $B_i$ may depend
on $i$.
Let $\sigma_i$ be the associated automorphism of $B_i$ in ${\mathcal B}$
($i=1,\ldots,m$).

\begin{prop} Under the above assumptions, $\sigma_i\sigma_j=\sigma_j\sigma_i$ on $\Pcl$
for any $i,j$.
\end{prop}

\begin{proof}
Take any $\la,\mu$ in $\Pcl^+$. Let $B'_i$ (resp. $B'_j$) be the perfect
crystal of level $\langle c,\la\rangle$ (resp. $\langle c,\mu\rangle$) which
belongs to the same coherent family as $B_i$ (resp. $B_j$). Let $B''_i$ be
the perfect crystal of level $\langle c,\la+\mu\rangle$ which belongs to the
same coherent family as $B_i$. {}From Theorems \ref{th:sum-of-B}, \ref{th:core}
and equation \eqref{eq:iso} one has the following isomorphisms
of $\Pcl$-weighted crystals
\begin{eqnarray*}
B(\la)\ot B(\mu)&\simeq&B(\sigma_i\la)\ot B'_i\ot B(\mu)
\simeq B(\sigma_i\la)\ot B(\sigma_i\mu)\ot B''_i\\
&\simeq&B(\sigma_i\la)\ot B(\sigma_j\sigma_i\mu)\ot B'_j\ot B''_i
\simeq B(\sigma_i\la)\ot B(\sigma_j\sigma_i\mu)\ot B''_i\ot B'_j.
\end{eqnarray*}
Going backward, we have $B(\la)\ot B(\mu)\simeq
B(\la)\ot B(\sigma_j^{-1}\sigma_i^{-1}\sigma_j\sigma_i\mu)$,
and therefore $\mu=\sigma_j^{-1}\sigma_i^{-1}\sigma_j\sigma_i\mu$.
{}From \eqref{eq:sigma-cond} we see that the statement is valid.
\end{proof}

Let $\tau$ be a Dynkin diagram automorphism, that is, a
permutation of the vertex set $I$ which preserves the Dynkin
diagram of $\geh$. By setting $\tau(\sum_i m_i\La_i) =\sum_i
m_i\La_{\tau(i)}$ one can extend it to a linear automorphism of
$\Pcl$. Since $\tau$ is a Dynkin automorphism it follows that
$\tau(\alpha_i)=\alpha_{\tau(i)}$ for all $i\in I$. It may happen
that $\tau$ can be lifted to an automorphism, denoted by
$\tau^\ast$, of a perfect crystal $B$ satisfying
\begin{equation} \label{eq:tau-on-B}
\tau^\ast(\et{i}b)=\et{\tau(i)}\tau^\ast(b)\mbox{ and }
\tau^\ast(\ft{i}b)=\ft{\tau(i)}\tau^\ast(b),
\end{equation}
for any $i$ and $b\in B$ such that $\et{i}b\neq0$ and $\ft{i}b\neq0$,
respectively.
If such a $\tau^\ast$ exists, let us call $\tau$ a {\em proper} automorphism
for $B$. Note that (\ref{eq:tau-on-B}) implies $\veps(\tau^\ast(b))=
\tau\veps(b),\vphi(\tau^\ast(b))=\tau\vphi(b)$ and hence $\tau^\ast$ is
unique for $B$.
It should be emphasized that the associated automorphism $\sigma$ of a
perfect crystal $B$ is an automorphism on $\Pcl$ and $\sigma$ can be a
proper automorphism for a different crystal $B'$ from $B$.
In what follows we simply write $\sigma b$ to mean $\sigma^\ast(b)$.

Recall $\sigma_i$ is the associated automorphism of $B_i$ in ${\mathcal B}$.

\begin{itemize}
\item[(A2)] For any pair $(i,j)$ ($1\le i<j\le m$) $\sigma_i$ is
induced by a Dynkin diagram automorphism (see Conjecture
\ref{conj:auto}(2)) and acts on $B_j$ as a proper automorphism.
Moreover, if $b_i\ot b_j\in(B_i\ot B_j)_{\min}$ and $b_i\ot b_j$
is mapped to $\bt_j\ot\bt_i$ under the isomorphism $B_i\ot
B_j\simarrow B_j\ot B_i$, then $\bt_j=\sigma_i b_j$.
\end{itemize}

\begin{remark} \label{rem:b-tilde}
The above $\bt_i$ is determined easily as $\bt_i=\vphi_{B_i}^{-1}
(\vphi_{B_i}(b_i)+\wt b_j)$ using (\ref{eq:ot-phi}) and
Lemma \ref{lem:tensor} (2).
\end{remark}

The following proposition is sometimes useful to check whether (A2) is
satisfied.

\begin{prop} \label{prop:combR}
Let $B_1,B_2$ be perfect crystals of level $l_1,l_2$ ($l_1\ge l_2$).
Suppose $B_1$ belongs to a certain coherent family of perfect crystals with the
associated automorphism $\sigma$ that is proper for $B_2$. We assume
\begin{equation} \label{eq:inj-assump}
\mbox{the map }\veps\times\vphi\;:\;B_2\longrightarrow(\Pcl^+)^2
\mbox{ is injective.}
\end{equation}
Then if $b_1\ot b_2\in(B_1\ot B_2)_{\min}$, $b_1\ot b_2$ is mapped to
$\sigma b_2\ot \bt_1$ under the isomorphism $B_1\ot B_2\simarrow B_2\ot B_1$.
\end{prop}

\begin{proof}
Let $b_1\ot b_2$ be mapped to $\bt_2\ot\bt_1$ under the isomorphism.
We are to show $\bt_2=\sigma b_2$. Set $\la=\sigma\veps(b_2)$,
$\mu=\vphi(b_2)$ and $B_0$ to be the perfect crystal
of level $l_1-\langle c,\la\rangle$ in the same coherent family as $B_1$.
Note that from Lemma \ref{lem:tensor} (2), we have
$\veps(b_1)=\sigma\vphi(b_1)\ge\sigma\veps(b_2)=\la$ and hence $l_1\ge
\langle c,\la\rangle$. Using Theorems \ref{th:sum-of-B}, \ref{th:core} and
Proposition \ref{prop:core} we have
\[
\begin{array}{ccccc}
B(\la)\ot B_1\ot B_2&\simeq&B_0\ot B(\sigma^{-1}\la)\ot B_2&\simeq&
B_0\ot(\oplus_\nu B(\nu)^{\oplus m_\nu})\\
u_\la\ot b_1\ot b_2&\mapsto&b_0\ot u_{\sigma^{-1}\la}\ot b_2&\mapsto&
b_0\ot u_\mu\\
\\
&\simeq&(\oplus_\nu B(\sigma\nu)^{\oplus m_\nu})\ot B_1\\
&\mapsto&u_{\sigma\mu}\ot b'_1,
\end{array}
\]
with some $b_0\in B_0,b'_1\in B_1$.
Now notice that when $B(\sigma^{-1}\la)\ot B_2\simeq
\oplus_\nu B(\nu)^{\oplus m_\nu}$, we have
$B(\la)\ot B_2\simeq\oplus_\nu B(\sigma\nu)^{\oplus m_\nu}$
because $\sigma$ is proper for $B_2$.
Since $b_1\ot b_2$ is mapped to $\bt_2\ot\bt_1$ under
$B_1\ot B_2\simarrow B_2\ot B_1$, by $B(\la)\ot B_2\ot B_1\simarrow
(\oplus_\nu B(\sigma\nu)^{\oplus m_\nu})\ot B_1$, $u_\la\ot\bt_2\ot\bt_1$
should be mapped to $u_{\sigma\mu}\ot\bt_1$. Therefore we have
\[
\la=\veps(\sigma b_2)=\veps(\bt_2),\quad
\sigma \mu=\vphi(\sigma b_2)=\vphi(\bt_2).
\]
{}From the assumption (\ref{eq:inj-assump}) we get $\bt_2=\sigma b_2$.
\end{proof}

\begin{example} \label{ex:KKM}
In \cite{KKM} a coherent family of perfect crystals ${\mathcal
B}=\{B^l\}_{l\ge1}$ is given for any nonexceptional affine Lie
algebra $\geh$. For these ${\mathcal B}$ (A1) is easily verified
using Lemma \ref{lem:auto}. Let $\sigma$ be the associated
automorphism. One can see this $\sigma$ is proper for any $B^l$,
since the automorphism on $B^l$ is determined uniquely by
requiring $\sigma b=b'$ if $\sigma\veps(b)=\veps(b')$ for $b,b'\in
B^l_{\min}$ and \eqref{eq:tau-on-B}. Using the explicit formula
for $\veps(b)$ and $\vphi(b)$ in \cite{KKM}, one can also check
that the condition \eqref{eq:inj-assump} holds. Thus by
Proposition \ref{prop:combR}, (A2) is also valid for this
${\mathcal B}$.
\end{example}

\begin{example} \label{ex:A}
Consider the set of perfect crystals
${\mathcal B}=\{B^{k,l}\}_{1\le k\le n,l\ge1}$ of type
$A^{(1)}_n$ given in \cite{KMN2}. We will show in the next section that
$\{B^{k,l}\}_{l\ge1}$ for any fixed $k$ forms a coherent family of perfect
crystals with the associated automorphism $\sigma_k$ satisfying
(\ref{eq:sigma-cond}) such that $\sigma_k(\La_i)=\La_{i-k\,\mathrm{mod}\,n+1}$.
Thus (A1) is valid for ${\mathcal B}$.
It is also known that there exists an automorphism $pr$ on any $B^{k,l}$ called
the promotion operator \cite{S}. $pr$ satisfies
\[
pr\circ\et{i-1}=\et{i}\circ pr\quad\mbox{and}\quad pr\circ\ft{i-1}=\ft{i}\circ pr
\]
for any $i$ where indices are taken modulo $n+1$. Therefore defining an
automorphism $\sigma_{k'}$ on $B^{k,l}$ by $\sigma_{k'}=pr^{-k'}$,
$\sigma_{k'}$ turns out to be a proper automorphism for $B^{k,l}$.
Furthermore, the second requirement in (A2) is already proven
in Theorem 7.3 of \cite{SS}. Hence (A2) is also valid for ${\mathcal B}$.
We note that the condition (\ref{eq:inj-assump}) is no longer true for
$B^{k,l}$ in general.
\end{example}

In what follows in this section, we assume that any set of perfect crystals we consider
satisfies (A1) and (A2).

\subsection{Lemmas}
We prepare several lemmas. First suppose $B_1$ and $B_2$ are perfect crystals of level
$l_1$ and $l_2$ such that $l_1\ge l_2$, and let $\sigma$ be the associated automorphism
of $B_1$.

\begin{lemma} \label{lem:linked2}
Let $b_1,b'_1\in B_1, b_2\in B_2$ and suppose $b_2\ot b_1\in(B_2\ot B_1)_{\min}$. Then there exists a sequence
$i_1,\ldots,i_l\in I$ such that
\[
b_2\ot b_1=\et{i_l}\cdots\et{i_1}(b_2\ot b'_1).
\]
\end{lemma}

\begin{proof}
Let $\la=\veps(b_2)$ and consider the connected component $\hat{B}$ of
$B(\la)\ot B_2\ot B_1$ containing $u_\la\ot b_2\ot b_1$. We have
$u_\la\ot b_2\ot b'_1\in\hat{B}$ and also $u_\la\ot b_2\ot b_1\in C(\hat{B})$
by Proposition \ref{prop:core} (2). Apply Proposition \ref{prop:core} (3).
\end{proof}

Set $r=\sigma^{-1}(0)$ and define a function $\gamma$ on $B_2$ by
\begin{equation} \label{eq:gamma}
\gamma(b_2)=(\La_r^\vee-\La_0^\vee|\wt b_2),\quad \La_i^\vee=(a_i/a_i^\vee)\La_i,
\end{equation}
where $a_i$ (resp. $a_i^\vee$) is the $i$-th Kac (resp. dual Kac) label and
$(\,|\,)$ is the invariant bilinear form on $P$ as in \cite{Kac}.
If $r=0$, $\gamma(b)=0$. Otherwise, it has the following recursive property
\begin{equation} \label{eq:def-gamma}
\gamma(\et{i}b_2)=\left\{
\begin{array}{ll}
\gamma(b_2)-1\quad&\mbox{if }i=0,\\
\gamma(b_2)+1&\mbox{if }i=r,\\
\gamma(b_2)&\mbox{otherwise}.
\end{array}\right.
\end{equation}

\begin{lemma} \label{lem:H=gamma}
If $b_1\ot b_2\in(B_1\ot B_2)_{\min}$, we have $H_{B_1B_2}(b_1\ot b_2)=\gamma(b_2)$
up to global additive constant.
\end{lemma}

\begin{proof}
Take $b_1\ot b_2,b'_1\ot b'_2\in(B_1\ot B_2)_{\min}$ and suppose $b_1\ot b_2\mapsto \bt_2\ot\bt_1,
b'_1\ot b'_2\mapsto\bt'_2\ot\bt'_1$ under the isomorphism $B_1\ot B_2\simarrow B_2\ot B_1$.
{}From (A2) we have $\bt_2=\sigma b_2,\bt'_2=\sigma b'_2$.
Take a sequence $i_1,\ldots,i_l\in I$ such that $\bt'_2=\et{i_l}\cdots\et{i_1}\bt_2$ and
$l$ is minimal. Then one can get a sequence $k_1,\ldots,k_n\in I$ such that
$\et{k_n}\cdots\et{k_1}(\bt_2\ot\bt_1)=\et{i_l}\cdots\et{i_1}\bt_2\ot b''_1$ for some $b''_1\in B_1$.
Using the previous lemma one can assume $b''_1=\bt'_1$.
Now suppose $\et{k_n}\cdots\et{k_1}(b_1\ot b_2)
=\et{i^{\prime\ast}_{l^{\prime\ast}}}\cdots\et{i^{\prime\ast}_1}b_1\ot\et{i'_{l'}}\cdots\et{i'_1}b_2$.
Due to the minimality of $l$,
\[
\sigma(\et{i'_{l'}}\cdots\et{i'_1}b_2)=\et{i_l}\cdots\et{i_1}\sigma b_2
\]
implies
\[
\sum_{a=1}^{l'}\alpha_{i'_a}-\sum_{a=1}^l\alpha_{\sigma^{-1}(i_a)}\in\Zn\delta.
\]
Similarly, going from $b'_1\ot b'_2$ to $b_1\ot b_2$ by applying $\et{i}$'s, one obtains
\begin{eqnarray*}
&&\et{k'_{n'}}\cdots\et{k'_1}(\bt'_2\ot\bt'_1)
=\et{j_m}\cdots\et{j_1}\bt'_2\ot\et{j^\ast_{m^\ast}}\cdots\et{j^\ast_1}\bt'_1,\\
&&\et{k'_{n'}}\cdots\et{k'_1}(b'_1\ot b'_2)
=\et{j^{\prime\ast}_{m^{\prime\ast}}}\cdots\et{j^{\prime\ast}_1}b'_1\ot\et{j'_{m'}}\cdots\et{j'_1}b'_2,\\
&&\sum_{a=1}^{m'}\alpha_{j'_a}-\sum_{a=1}^m\alpha_{\sigma^{-1}(j_a)}\in\Zn\delta.
\end{eqnarray*}
Calculating the difference of the energy function during this process from Lemma \ref{lem:diff-H},
one gets
\[
0=\sharp0[i_1,\ldots,i_l,j_1\ldots,j_m]-\sharp0[i'_1,\ldots,i'_{l'},j'_1\ldots,j'_{m'}].
\]
Therefore we have $l'=l,m'=m$ and
\[
[i'_1,\ldots,i'_{l'}]=[\sigma^{-1}(i_1),\ldots,\sigma^{-1}(i_l)]\quad\mbox{as multiset}.
\]
Thus noting $b'_2=\et{\sigma^{-1}(i_l)}\cdots\et{\sigma^{-1}(i_1)}b_2$, one gets
\begin{eqnarray*}
H(b'_1\ot b'_2)-H(b_1\ot b_2)
&=&\sharp0[i_1,\ldots,i_l]-\sharp0[\sigma^{-1}(i_1),\ldots,\sigma^{-1}(i_l)]\\
&=&\sharp{r}[\sigma^{-1}(i_1),\ldots,\sigma^{-1}(i_l)]-\sharp0[\sigma^{-1}(i_1),\ldots,\sigma^{-1}(i_l)]\\
&=&\gamma(b'_2)-\gamma(b_2)
\end{eqnarray*}
as desired.
\end{proof}

Let ${\mathcal B}=\{B_1,B_2,\ldots,B_m\}$ be a set of perfect
crystals satisfying (A1), (A2) and $\lev B_1\ge\lev
B_2\ge\cdots\ge\lev B_m$.

\begin{lemma} \label{lem:propagation}
Suppose $b_1\ot\cdots\ot b_m\in(B_1\ot\cdots\ot B_m)_{\min}$. Then
$b_1^{\langle j-1\rangle}\ot b_j\in(B_1\ot B_j)_{\min}$ for any $2\le j\le m$.
Here $b_1^{\langle j\rangle}$ is defined in \eqref{eq:b^<j>}.
\end{lemma}

\begin{proof}
Let $b_1\ot b_2\ot\cdots\ot b_{j-1}$ be mapped to $\bt_2\ot\cdots\ot\bt_{j-1}\ot
b^{\langle j-1\rangle}_1$ by $B_1\ot B_2\ot\cdots\ot B_{j-1}\simarrow
B_2\ot\cdots\ot B_{j-1}\ot B_1$. From the fact that
$\bt_2\ot\cdots\ot\bt_{j-1}\ot b^{\langle j-1\rangle}_1\ot b_j\cdots\ot b_m$
is minimal and using Lemma \ref{lem:tensor}, (\ref{eq:ot-eps}),
(\ref{eq:ot-phi}), we have
\[
\vphi(b^{\langle j-1\rangle}_1)=\vphi(\bt_2\ot\cdots\ot\bt_{j-1}\ot
b^{\langle j-1\rangle}_1)\ge\veps(b_j\ot\cdots\ot b_m)\ge\veps(b_j)
\]
and know $b^{\langle j-1\rangle}_1$ is minimal. Hence
$b^{\langle j-1\rangle}_1\ot b_j\in(B_1\ot B_j)_{\min}$ again by
Lemma \ref{lem:tensor} (2).
\end{proof}

Let $\sigma$ be the associated automorphism of $B_1$ and
set $\Bd=B_2\ot\cdots\ot B_m$. For an element $\bd=b_2\ot\cdots\ot b_m\in
\Bd$ define $\sigma\bd$ by
\[
\sigma\bd=\sigma b_2\ot\cdots\ot\sigma b_m.
\]
Note that this definition allows (\ref{eq:tau-on-B}) to hold on the tensor
product crystal. Let $\gamma_{B_j}$ be the function $\gamma$ (\ref{eq:gamma})
on $B_j$. The following lemma is now immediate from Lemma \ref{lem:propagation},
(A2), Remark \ref{rem:b-tilde}, Proposition \ref{prop:HofTensor} (2) and
Lemma \ref{lem:H=gamma}.

\begin{lemma} \label{lem:B1-Bd}
Let $b_1\ot\bd\in(B_1\ot\Bd)_{\min}$. Then

\begin{itemize}
\item[(1)] $b_1\ot\bd$ is mapped to $\sigma\bd\ot\bt_1$ under the isomorphism
           $B_1\ot\Bd\simarrow\Bd\ot B_1$. Here
           $\bt_1=\vphi_{B_1}^{-1}(\vphi_{B_1}(b_1)+\wt\bd)$.

\item[(2)] $H_{B_1\Bd}(b_1\ot\bd)=\gamma_{\Bd}(\bd)
            :=\sum_{j=2}^m\gamma_{B_j}(b_j)$ up to global additive constant.
\end{itemize}
\end{lemma}

Next we consider the tensor product $B_1^{\ot L}\ot\cdots\ot B_m^{\ot L}$.
Set $B_1^\star=B_1^{\ot L},B^\ddagger=B_2^{\ot L}\ot\cdots\ot B_m^{\ot L}$.
Here is the last lemma.

\begin{lemma} \label{lem:last}
Let $L$ be a multiple of the order of $\sigma$.
Suppose $b_1\in(B_1^\star)_{\min}$. Then the map $B^\ddagger\longrightarrow
\Z$ given by $b^\ddagger\mapsto H_{B_1^\star B^\ddagger}(b_1\ot b^\ddagger)$
is constant on the set of elements $b^\ddagger$ such that $b_1\ot b^\ddagger$
is minimal.
\end{lemma}

\begin{proof}
Write $b_1=b^1\ot\cdots\ot b^L,b^i\in B_1$. By Proposition
\ref{prop:HofTensor} (1) and Lemma \ref{lem:B1-Bd} we have
\[
H_{B_1^\star B^\ddagger}(b_1\ot b^\ddagger)
=\sum_{i=1}^L H_{B_1 B^\ddagger}(b^{L+1-i}\ot\sigma^{i-1} b^\ddagger)
=c+\sum_{i=1}^L \gamma_{B^\ddagger}(\sigma^{i-1}b^\ddagger)
\]
where $c$ is a constant. However for any weight $\la$, one has
\[
(\La^\vee_r\mid\sum_{i=1}^L\sigma^{i-1}\la)
=(\sigma\La^\vee_r\mid\sum_{i=1}^L\sigma^i\la)
=(\La^\vee_0\mid\sum_{i=1}^L\sigma^{i-1}\la).
\]
Hence $\sum_{i=1}^L\gamma_{B^\ddagger}(\sigma^{i-1}b^\ddagger)=0$ by
\eqref{eq:gamma}.
\end{proof}

\subsection{Main theorem}
Let ${\mathcal B}=\{B_1,B_2,\ldots,B_m\}$ be a set of perfect crystals
satisfying (A1) and (A2). Set $l_i=\lev B_i$ and assume
$l_1\ge l_2\ge\cdots\ge l_m\ge l_{m+1}:=0$.
Let $\sigma_i$ be the associated automorphism of $B_i$ and $B'_i$ be a
perfect crystal of level $l_i-l_{i+1}$
which belongs to the same coherent family as $B_i$. Take
$\la_i\in(\Pcl^+)_{l_i-l_{i+1}}$ for all $i=1,2,\ldots,m$.
Then we have the isomorphisms of $\Pcl$-weighted crystals
\begin{eqnarray}
&&B(\la_1)\ot B(\la_2)\ot\cdots\ot B(\la_m) \label{eq:main iso}\\
&\simeq&(B(\sigma_1\la_1)\ot B'_1)\ot(B(\sigma_2\la_2)\ot B'_2)\ot\cdots\ot(B(\sigma_m\la_m)\ot B'_m) \nonumber\\
&\simeq&B(\sigma_1\la_1)\ot B(\sigma_1\sigma_2\la_2)\ot\cdots\ot B(\sigma_1\sigma_2\cdots\sigma_m\la_m)
\ot(B_1\ot B_2\ot\cdots\ot B_m) \nonumber
\end{eqnarray}
by applying Theorems \ref{th:sum-of-B} and \ref{th:core} successively.

\begin{prop} \label{prop:shippo}
Under the isomorphism \eqref{eq:main iso} $u_{\la_1}\ot u_{\la_2}\ot\cdots\ot u_{\la_m}$ is sent to
$u_{\sigma_1\la_1}\ot u_{\sigma_1\sigma_2\la_2}\ot\cdots\ot u_{\sigma_1\sigma_2\cdots\sigma_m\la_m}
\ot(\bb^{(1)}\ot\bb^{(2)}\ot\cdots\ot\bb^{(m)})$ ($\bb^{(i)}\in B_i$) where $\bb^{(i)}$ is given by
\[
\bb^{(i)}=\veps_{B_i}^{-1}(\sigma_i\la_i+\sigma_i\sigma_{i+1}\la_{i+1}+\cdots+\sigma_i\sigma_{i+1}\cdots\sigma_m\la_m)
\]
for $i=1,2,\ldots,m$.
\end{prop}

\begin{proof}
We consider the $m=2$ case first. Under \eqref{eq:main iso} for $m=2$ we have
\[
u_{\la_1}\ot u_{\la_2}\mapsto(u_{\sigma_1\la_1}\ot\bb^{\prime(1)})\ot(u_{\sigma_2\la_2}\ot\bb^{\prime(2)})
\mapsto u_{\sigma_1\la_1}\ot u_{\sigma_1\sigma_2\la_2}\ot(\bb^{(1)}\ot\bb^{(2)}).
\]
To determine $\bb^{\prime(i)}$, calculate $\veps(u_{\sigma_i\la_i}\ot\bb^{\prime(i)})$, which must be $0$, using \eqref{eq:ot-eps}.
We know $\bb^{(2)}=\bb^{\prime(2)}$. To obtain $\bb^{(1)}$, use the fact that $\veps$ is unchanged under
$B'_1\ot B(\sigma_2\la_2)\simarrow B(\sigma_1\sigma_2\la_2)\ot B_1:\bb^{\prime(1)}\ot u_{\sigma_2\la_2}\mapsto
u_{\sigma_1\sigma_2\la_2}\ot\bb^{(1)}$.

The cases when $m>2$ are similar.
\end{proof}

Iterating the isomorphism \eqref{eq:main iso} we are led to consider a set of paths based on the finite crystal
$B=B_1\ot B_2\ot\cdots\ot B_m$. In view of Proposition \ref{prop:shippo} we define a distinguished path
$\pb^{(\la_1,\ldots,\la_m)}$ associated with a set of weights $\la_1,\ldots,\la_m$ ($\la_i\in(\Pcl^+)_{l_i-l_{i+1}}$) by
\begin{eqnarray}
\pb^{(\la_1,\ldots,\la_m)}(j)&=&\bb^{(1)}_j\ot\cdots\ot\bb^{(m)}_j\quad(j\ge1,\bb^{(i)}_j\in B_i\text{ for }1\le i\le m),
\label{eq:ref path}\\
\bb^{(i)}_j&=&\veps_{B_i}^{-1}(\sigma_1^{j-1}\cdots\sigma_{i-1}^{j-1}
\sum_{k=i}^m\sigma_i^j\cdots\sigma_k^j\la_k). \nonumber
\end{eqnarray}
{}From Proposition \ref{prop:shippo} the highest weight element of
$B(\la_1)\ot\cdots\ot B(\la_m)$ is sent to
$u_{\sigma_1^L\la_1}\ot\cdots\ot u_{\sigma_1^L\cdots\sigma_m^L\la_m}\ot
\pb^{(\la_1,\ldots,\la_m)}(L)\ot\cdots\ot\pb^{(\la_1,\ldots,\la_m)}(1)$
under the $L$ fold iteration of \eqref{eq:main iso}, which implies that
$\pb^{(\la_1,\ldots,\la_m)}(L)\ot\cdots\ot\pb^{(\la_1,\ldots,\la_m)}(1)$ is
minimal in $B^{\ot L}$. This fact shows

\begin{prop}
$\pb^{(\la_1,\ldots,\la_m)}$ is a reference path in $B_1\ot\cdots\ot B_m$.
\end{prop}

We now have a set of paths $\P(\pb^{(\la_1,\ldots,\la_m)},B_1\ot\cdots\ot B_m)$.
Our main theorem is the following.

\begin{theorem} \label{th:main}
Let $B_i$ be a perfect crystal of level $l_i$ with the associated automorphism
$\sigma_i$ for $i=1,\dots,m$. Suppose $l_1\ge l_2\ge\cdots\ge l_m$.
For $\la_i\in(\Pcl^+)_{l_i-l_{i+1}}$ ($i=1,\dots,m;l_{m+1}=0$) let
$\pb^{(\la_1,\dots,\la_m)}$ be a reference path in $B_1\ot\cdots\ot B_m$ defined in \eqref{eq:ref path}.
Then we have an isomorphism of $P$-weighted crystals
\[
\P(\pb^{(\la_1,\dots,\la_m)},B_1\ot\cdots\ot B_m)\simeq
B(\la_1)\ot\cdots\ot B(\la_m).
\]
\end{theorem}

\begin{proof}
We have the isomorphism of $\Pcl$-weighted crystals by iterating \eqref{eq:main iso}.
\[
B(\la_1)\ot\cdots\ot B(\la_m)
\simeq B(\sigma_1^L\la_1)\ot\cdots\ot B(\sigma_1^L\cdots\sigma_m^L\la_m)\ot(B_1\ot\cdots\ot B_m)^{\ot L}
\]
For an element $v_1\ot\cdots\ot v_m\in B(\la_1)\ot\cdots\ot B(\la_m)$ take sufficiently large $L$ such that
under the above isomorphism
\[
v_1\ot\cdots\ot v_m\mapsto u\ot p(L)\ot\cdots\ot p(1),
\]
where $u=u_{\sigma_1^L\la_1}\ot\cdots\ot u_{\sigma_1^L\cdots\sigma_m^L\la_m}$ and
$p(j)\in B_1\ot\cdots\ot B_m$.
We also assume that $L(>0)$ is a multiple of the order of $\sigma_i$ for all
$i$ and $p(L)=\pb(L)$ where
$\pb(L)$ is the $L$-th component of the reference path $\pb=\pb^{(\la_1,\ldots,\la_m)}$.
To prove the theorem it suffices to show
\begin{equation} \label{eq:to show}
-\sum_{i=1}^m\langle d,\wt v_i\rangle=\sum_{j=1}^{L-1}j(H(p(j+1)\ot p(j))-H(\pb(j+1)\ot\pb(j))),
\end{equation}
where $d$ is the scaling element of $\geh$ as in \cite{Kac} and $H$ is the energy function on
$B\ot B$ ($B=B_1\ot\cdots\ot B_m$).
See \eqref{eq:def E}, \eqref{eq:bold wt} and recall $\langle d,\delta\rangle=a_0$.

Consider the isomorphism of crystals
\[
B^{\ot L}\simarrow B_1^{\ot L}\ot\cdots\ot B_m^{\ot L},
\]
where $B=B_1\ot\cdots\ot B_m$. Suppose $p(L)\ot\cdots\ot p(1)\mapsto b_1\ot\cdots\ot b_m,
\pb(L)\ot\cdots\ot\pb(1)\mapsto\bb_1\ot\cdots\ot\bb_m$ ($b_i,\bb_i\in B_i^{\ot L}$) under the isomorphism.
By \eqref{eq:hom Denergy} and Proposition \ref{prop:D intrinsic} we have
\[
r.h.s.\text{ of \eqref{eq:to show}}=D(b_1\ot\cdots\ot b_m)-D(\bb_1\ot\cdots\ot\bb_m).
\]

We wish to show \eqref{eq:to show} by induction on $m$. We first observe that it suffices to check
\eqref{eq:to show} when $v_1\ot\cdots\ot v_m$ is a highest weight element. If not, we may certainly
apply $\et{i}$ for some $i$ with the change $-\delta_{i0}$ on both sides of \eqref{eq:to show} from
\eqref{eq:hom Denergy} and Lemma \ref{lem:D diff}. Thus we may assume that $v_1\ot\cdots\ot v_m$ 
is a highest weight element.
This implies $v_1=u_{\la_1}$ and that $b_1\ot\cdots\ot b_m$ is a minimal element of
$B_1^{\ot L}\ot\cdots\ot B_m^{\ot L}$, and therefore that $b_1=\bb_1$. Write $B_1^\star=B_1^{\ot L},
B^\ddagger=B_2^{\ot L}\ot\cdots\ot B_m^{\ot L},b^\ddagger=b_2\ot\cdots\ot b_m$ and
$\bb^\ddagger=\bb_2\ot\cdots\ot\bb_m$. We have
\begin{eqnarray*}
D(b_1\ot b^\ddagger)-D(\bb_1\ot\bb^\ddagger)&=&
D(\bb_1\ot b^\ddagger)-D(\bb_1\ot\bb^\ddagger)\\
&=&H_{B_1^\star B^\ddagger}(\bb_1\ot b^\ddagger)-H_{B_1^\star B^\ddagger}(\bb_1\ot\bb^\ddagger)+D(b^\ddagger)-D(\bb^\ddagger)\\
&=&D(b^\ddagger)-D(\bb^\ddagger)
\end{eqnarray*}
by \eqref{eq:func} for two tensor factors,
(A2), the fact that
the order of $\sigma_1$ divides $L$, and Lemma \ref{lem:last}.

Noting $\langle d,\wt u_{\la_1}\rangle=0$ we have reduced
\eqref{eq:to show} to the case that the first tensor factor is missing,
and the theorem is proved by induction.
\end{proof}

\section{Type $A$ perfect crystals} \label{sec:A}

In this section we review the perfect crystal $B^{k,l}$ of type $A^{(1)}_n$,
given in \cite{KMN2}, and show for fixed $k$ ($1\le k\le n$) that
$\{B^{k,l}\}_{l\ge1}$ forms a coherent family of perfect crystals.

\subsection{Crystal $B^{k,l}$}
Let $B^{k,l}$ ($1\le k\le n,l\ge1$) be the perfect crystal of level $l$ of
type $A^{(1)}_n$ given in \cite{KMN2}. As a $U_q(A_n)$-crystal it coincides
with the crystal base $B(l\ol{\La}_k)$ of the highest weight $U_q(A_n)$-module
with highest weight $l\ol{\La}_k$, where $\ol{\La}_k$ is the $k$-th fundamental
weight of $A_n$. Therefore $B^{k,l}$ can be identified with the set of
semistandard tableaux of $k\times l$ rectangle shape over the alphabet
$\{1,2,\ldots,n+1\}$ as described in \cite{KN}. Thus with each $b\in B^{k,l}$
one associates a table $(m_{jj'})_{1\le j\le k,1\le j'\le l}$ such that
$m_{jj'}\in\{1,2,\ldots,n+1\}, m_{jj'}\le m_{j,j'+1}$ and $m_{jj'}<m_{j+1,j'}$.
Set $k'=n+1-k$. For our purpose we associate another table $x=x(b)$ with $b$, namely
\begin{equation}
\label{eq:x_ji}
\begin{split}
x&=(x_{ji})_{1\le j\le k,j\le i\le j+k'},\\
x_{ji}&=\sharp\{j'\mid m_{jj'}=i\}.
\end{split}
\end{equation}
Notice that $x_{ji}$ defined by (\ref{eq:x_ji}) satisfies
\[
x_{ji}=0\quad\text{unless}\quad j\le i\le j+k'
\]
due to the semistandardness of the tableau. In view of this we set $x_{ji}=0$ when
$1\le i<j$ or $j+k'<i\le n+1$. We also set $x_{0i}=x_{k+1,i}=0$
($1\le i\le n+1$) for convenience.

It is shown in \cite{KMN2} that $B^{k,l}$ is perfect of level $l$.
We give below its associated automorphism $\sigma=\sigma_k$, which depends
only on $k$, and the minimal elements. Set
$\sigma_k(\La_i)=\La_{i-k\,\mathrm{mod}\,n+1}$. $\sigma_k$ on $\Pcl$ is
given by extending it $\Z$-linearly. Let $\la=\sum_{i=0}^n\la_i\La_i$ be
in $\Pcll$, that is, $\la_0,\la_1,\ldots,\la_n\in\Zn,\sum_{i=0}^n\la_i=l$.
The table $x(b)$ of the minimal element $b$ such that $\veps(b)=\la$ is
given by
\begin{equation}
\label{eq:min elem}
\begin{split}
x_{jj}&=\la_0+\sum_{\alpha=j}^{k-1}\la_{\alpha+k'},\\
x_{ji}&=\la_{i-j}\quad(j<i<j+k'),\\
x_{j,j+k'}&=\sum_{\alpha=0}^{j-1}\la_{\alpha+k'}
\end{split}
\end{equation}
for $1\le j\le k$.

\subsection{Actions of $\et{a},\ft{a}$ ($a\neq0$)}
We give the actions of $\et{a},\ft{a}$ ($a\neq0$) in terms of the coordinate $x(b)=(x_{ji})$.
Set $\beta=\max(0,a-k'),\gamma=\min(k,a)$. For fixed $x=x(b)$ define
\begin{equation} \label{eq:def Gamma}
\Gamma(c)=\sum_{j=\beta+1}^{c-1}(x_{ja}-x_{j+1,a+1})
\end{equation}
for $\beta+1\le c\le\gamma$. Let $\Gamma_{\min}$ be the minimum of $\Gamma(c)$.
Set
\begin{eqnarray*}
c_0&=&\min\{c\mid \beta+1\le c\le\gamma,\Gamma(c)=\Gamma_{\min}\},\\
c_1&=&\max\{c\mid \beta+1\le c\le\gamma,\Gamma(c)=\Gamma_{\min}\}.
\end{eqnarray*}
Then the values $\veps_a(b),\vphi_a(b)$ are given by
\begin{eqnarray}
\veps_a(b)&=&\sum_{j=\beta}^{c_0-1}(x_{j+1,a+1}-x_{ja}), \label{eq:Bkl eps}\\
\vphi_a(b)&=&\sum_{j=c_1}^\gamma(x_{ja}-x_{j+1,a+1}). \label{eq:Bkl phi}
\end{eqnarray}
If $\veps_a(b)>0$, writing $x'=x(\et{a}b)$ we have
\begin{equation} \label{eq:Bkl e}
x'_{ji}=x_{ji}-\delta_{jc_0}\delta_{i,a+1}+\delta_{jc_0}\delta_{ia}.
\end{equation}
Also if $\vphi_a(b)>0$, writing $x'=x(\ft{a}b)$ we have
\begin{equation} \label{eq:Bkl f}
x'_{ji}=x_{ji}-\delta_{jc_1}\delta_{ia}+\delta_{jc_1}\delta_{i,a+1}.
\end{equation}
These formulas are obtained by interpreting the rule
on the tableau given in \cite{KN} in terms of the coordinate $x(b)$.

\subsection{Actions of $\et{0},\ft{0}$} \label{ssec:KMN}
The following method (which we shall call the KKMMNN method) for
computing $\et{0}$ was stated in \cite[Prop. 6.3.11]{KMN2} but the
proof was not included. Proofs of this method and \eqref{eq:Bkl eps0} and
\eqref{eq:Bkl phi0} are given in Appendix \ref{sec:e0proof}.

The set $C$ of all sequences $1=c_0<c_1<\dotsm<c_{k-1}<c_k=n+1$ is
a poset under the relation $c \subseteq c'$ defined by $c_i\le
c'_i$ for all $0\le i\le k$. $C$ contains a unique
$\subseteq$-minimum element $c_{\min}=(1,2,\dotsc,k,n+1)\in C$ and
a unique $\subseteq$-maximum element
$c_{\max}=(1,k'+1,\dotsc,n+1)$ where $k'=n+1-k$.

Fix $x=x(b)$. Let $\Delta_x=\Delta:C\rightarrow\Zn$ be given by
\begin{equation}\label{eq:delta def}
\Delta(c) = \sum_{j=1}^k \sum_{c_{j-1}<i<c_j} x_{ji}.
\end{equation}
Consider the nonempty subposet $C^x_{\min}=C_{\min}$ of $C$
consisting of elements $c$ such that $\Delta(c)$ is minimal. The
poset $C_{\min}$ has a unique $\subseteq$-minimum element. For
suppose not. Let $c,c'\in C_{\min}$ be $\subseteq$-minima which
are incomparable. Define $c\wedge c'\in C$ and $c\vee c'\in C$ be
defined by
\[
\begin{split}
  (c\wedge c')_i &= \min(c_i,c'_i) \\
  (c\vee c')_i &= \max(c_i,c'_i)
\end{split}
\]
for all $i$. It can be verified that
\begin{equation} \label{eq:diffs}
  \Delta(c\wedge c') - \Delta(c) = \Delta(c') - \Delta(c\vee c').
\end{equation}
Since $c$ and $c'$ are incomparable, $c\wedge c'$ is properly
smaller than $c$. Thus the left hand side of \eqref{eq:diffs} is
positive. But the positivity of the right hand side contradicts
$c'\in C_{\min}$.

So for a given $x=x(b)$, there is a unique
$c\in C$ such that
\begin{align}\label{eq:mas}
\Delta(c) &\le \Delta(m) \qquad \text{if $m \supseteq c$} \\
\label{eq:mass1} \Delta(c) &< \Delta(m) \qquad \text{if
$m\not\supseteq c$}
\end{align}
for all $m\in C$.
The formula of \cite{KMN2} states that
\begin{enumerate}
\item $\et{0}b=0$ if and only if $c=c_{\min}$
and $x_{kk}=0$.
\item If $\et{0}b\neq0$, then writing $x'=x(\et{0}b)$, one has
\begin{equation} \label{eq:Bkl e0}
  x'_{ji} = x_{ji} - \delta_{i,c_{j-1}} + \delta_{i,c_j}
  \qquad\text{for all $1\le j\le k$ and $1\le i\le n+1$.}
\end{equation}
\end{enumerate}

The value $\veps_0(b)$ is given by
\begin{equation} \label{eq:Bkl eps0}
\veps_0(b)=l-x_{k,n+1}-\Delta(c).
\end{equation}

The formula for $\ft{0}$ is similar. It can be shown as above,
that for a given $x=x(b)$, there is a unique $\subseteq$-maximum
element $c$ in $C_{\min}$, that is, a unique $c\in C$ such that
\begin{align*}
\Delta(c) &\le \Delta(m) \qquad \text{if $m \subseteq c$} \\
\Delta(c) &< \Delta(m) \qquad \text{if
$m\not\subseteq c$}
\end{align*}
for all $m\in C$. Then $\ft{0}b$ is defined as follows.
\begin{enumerate}
\item $\ft{0}b=0$ if and only if $c=c_{\max}$
and $x_{1,k'+1}=0$.
\item If $\ft{0}b\neq0$, then writing $x'=x(\ft{0}b)$, one has
\begin{equation} \label{eq:Bkl f0}
  x'_{ji} = x_{ji} - \delta_{i,c_j} + \delta_{i,c_{j-1}}
  \qquad\text{for all $1\le j\le k$ and $1\le i\le n+1$.}
\end{equation}
\end{enumerate}
One also has
\begin{equation} \label{eq:Bkl phi0}
\vphi_0(b)=l-x_{11}-\Delta(c).
\end{equation}

\subsection{Coherent family}
We define the crystal $B^{k,\infty}$ by
\[
B^{k,\infty}=\{(\nu_{ji})_{1\le j\le k,j\le i\le j+k'}\mid
\nu_{ji}\in\Z,\sum_{i=j}^{j+k'}\nu_{ji}=0\text{ for any }j\}
\]
with $b_\infty=(\nu_{ji}=0\text{ for any }j,i)$. We again set
$\nu_{0i}=\nu_{k+1,i}=0$ ($1\le i\le n+1$) for convenience.
For $b=(\nu_{ji})\in B^{k,\infty}$, $\veps_a(b),\vphi_a(b),
\et{a}b,\ft{a}b$ for $a=1,2,\ldots,n$ are defined by the formulas
\eqref{eq:Bkl eps},\eqref{eq:Bkl phi},\eqref{eq:Bkl e},\eqref{eq:Bkl f}
with $x_{ji}$ replaced by $\nu_{ji}$, but without the condition $\veps_a(b)>0$
before \eqref{eq:Bkl e} and $\vphi_a(b)>0$ before \eqref{eq:Bkl f}. For $\et{0}$ and $\ft{0}$ use
\eqref{eq:Bkl e0} and \eqref{eq:Bkl f0}, but we warn that the formulas
for $\veps_0(b)$ and $\vphi_0(b)$ are modified as
\begin{equation} \label{eq:Binf 0}
\veps_0(b)=-\nu_{k,n+1}-\Delta(c),\quad\vphi_0(b)=-\nu_{11}-\Delta(c).
\end{equation}
Note that $\et{a}b,\ft{a}b\neq0$ for any $a=0,1,\ldots,n$ and
$b\in B^{k,\infty}$. Clearly we have $\veps_a(b_\infty)=\vphi_a(b_\infty)=0$
for any $a=0,1,\ldots,n$.

We are to show that $B^{k,\infty}$ is a limit of $\{B^{k,l}\}_{l\ge1}$
in the sense of section \ref{subsec:perfect}. For $b_0\in\Bmin^{k,l}$
let $x(b_0)=(\xi_{ji})$ and $\veps(b_0)=\la=\sum_{\alpha=0}^n\la_\alpha
\La_\alpha$. See \eqref{eq:min elem} for the explicit value of $\xi_{ji}$.
Then $\vphi(b_0)=\sigma_k^{-1}\la
=\sum_{\alpha=0}^{k-1}\la_{\alpha+k'}\La_\alpha
+\sum_{\alpha=k}^n\la_{\alpha-k}\La_\alpha$.
For $b\in B^{k,l}$ such that $x(b)=(x_{ji})$ we define the map
\[
f_{(l,b_0)} : T_\la\ot B^{k,l}\ot T_{-\sigma_k^{-1}\la}\longrightarrow
B^{k,\infty}
\]
by
\[
f_{(l,b_0)}(t_\la\ot b\ot t_{-\sigma_k^{-1}\la})=b'=(\nu_{ji}),
\]
where
\[
\nu_{ji}=x_{ji}-\xi_{ji}\quad\text{for }1\le j\le k,j\le i\le j+k'.
\]
Note that $f_{(l,b_0)}(t_\la\ot b_0\ot t_{-\sigma_k^{-1}\la})=b_\infty$.

Let us show that $f_{(l,b_0)}$ is a morphism of crystals. For this purpose
we prepare some properties of the functions $\Gamma$ \eqref{eq:def Gamma}
and $\Delta$ \eqref{eq:delta def}. Writing their $x$-dependence as
$\Gamma_x$ and $\Delta_x$, one can check
\begin{eqnarray*}
&&\Gamma_x(c)=\Gamma_{x-\xi}(c),\\
&&\Delta_x(c)-\Delta_{x-\xi}(c)=\sum_{\alpha=1}^{k'-1}\la_\alpha.
\end{eqnarray*}
Here $x-\xi$ means $(x_{ji})$ is replaced by $(x_{ji}-\xi_{ji})$.
{}From these formulas it is straightforward to see that
$f_{(l,b_0)}$ satisfies \eqref{eq:morph ef}. By \eqref{eq:Bkl
eps},\eqref{eq:Bkl phi},\eqref{eq:Bkl eps0},\eqref{eq:Bkl
phi0},\eqref{eq:Binf 0} and Example \ref{ex:T}, we also see that
\begin{equation*}
\begin{split}
\veps_a(t_\la\ot b\ot t_{-\sigma_k^{-1}\la})&=
\veps_a(b)-\langle h_a,\la\rangle=\veps_a(b'),\\
\vphi_a(t_\la\ot b\ot t_{-\sigma_k^{-1}\la})&=
\vphi_a(b)+\langle h_a,-\sigma_k^{-1}\la\rangle=\vphi_a(b')
\end{split}
\end{equation*}
for $a=0,1,\ldots,n$.

Hence for any $b_0\in\Bmin^{k,l}$, $f_{(l,b_0)}:
T_\la\ot B^{k,l}\ot T_{-\sigma_k^{-1}\la}\longrightarrow B^{k,\infty}$
is a morphism of crystals. By definition, it is clear that
$f_{(l,b_0)}$ is an embedding and that
\[
B^{k,\infty}=\bigcup_{(l,b_0)}\mbox{Im }f_{(l,b_0)}.
\]
Therefore $B^{k,\infty}$ is the limit of the coherent family of
perfect crystals $\{B^{k,l}\}_{l\ge1}$.

\appendix

\section{Formulas for $\et{0}$ and $\ft{0}$}
\label{sec:e0proof}

In this appendix we prove that two ways of
computing $\et{0}b$ for the perfect crystal $B^{k,l}$ of type
$A^{(1)}_n$, one given by \cite{KMN2} and the other by \cite{S},
are equivalent. The formulas for $\veps_0,\vphi_0$ and $\ft{0}$ are
also proven.

\subsection{Promotion method}
\label{ssec:jeu} In this subsection the method of \cite{S} to compute
$\et{0}$ is reviewed and suitably reformulated.

Let $b\in B^{k,\ell}$. The following algorithm to compute
$\et{0}b$ is easily seen to be equivalent to that in \cite{S}. We
shall use the notation $x$ interchangeably with $b$ when $x=x(b)$.
Let $x_j$ denote the $j$-th row. $x$ can also be identified with
its row word $x_kx_{k-1}\cdots x_1$.

\begin{enumerate}
\item Remove the letters $n+1$ from $x$, forming the subtableau
$x|_{[n]}$. Since the shape of $x$ is the $k\times \ell$
rectangle, all of these letters $n+1$ lie in the last row (the
$k$-th). Let $m=x_{k,n+1}$ be the number of such letters $n+1$.
\item Slide $x|_{[n]}$ to antinormal shape, obtaining the skew
tableau $z$ of shape $(\ell^k)/(m)$. Let $z_j$ be the $j$-th row
of $z$.
\item $\et{0}x$ is defined if and only if $z_1$ contains a letter
$1$.
\item Suppose this holds. Let $z'$ be obtained from $z$ by removing the
leftmost letter $1$ in $z_1$.
\item Slide the skew tableau $z'$ to normal shape. It has
shape $(\ell^{k-1},\ell-m-1)$.
\item Append $m+1$ letters $n+1$ to the last row; the result is $\et{0}x$.
\end{enumerate}
Note that in passing from $b$ to $\et{0}b$, a letter $1$ has been
removed and a letter $n+1$ added, with the other letters moving
around.

It is useful to perform the sliding algorithms in a particular
way.
Define the row words (weakly increasing words) $y_j$ for $1\le
j\le k$ as follows. Let $y_k=(x_k)|_{[n]}$, that is, let $y_k$ be
obtained by removing the letters $n+1$ from the last row $x_k$ of
$x$. Then $x|_{[n]}=y_k x_{k-1} \dotsm x_2 x_1$. Inductively (for
$j$ decreasing from $k-1$ down to $1$) let $y_j$ be the unique row
word of length $\ell-m$ and $z_{j+1}$ the unique row word of
length $\ell$, such that
\begin{equation} \label{eq:twojeu}
  y_{j+1} x_j \equiv z_{j+1} y_j.
\end{equation}
Also define $z_1=y_1$. Then the tableau $z$ defined above, is
given by $z=z_k \dotsm z_2 z_1$. The entries of $y_j$ are given
explicitly by the following relations.
\begin{equation} \label{eq:y}
\begin{aligned}
y_{ki}&=x_{ki} &\text{for $1\le i<n+1$}\\
y_{k n+1}&=0&\\
y_{ji}&=\min(x_{ji},\sum_{\alpha>i} (y_{j+1
\alpha}-y_{j\alpha})).&\qquad\text{for $1\le j\le k-1$.}
\end{aligned}
\end{equation}
The statements on $y_k$ are immediate. To justify the equation for
$y_j$ with $1\le j\le k-1$, note that \eqref{eq:twojeu} can be
computed by a two-row jeu de taquin. Consider the two-row tableau
$y_{j+1} x_j$. We imagine that $m$ letters are moving out of $x_j$
to the other row, to create $y_j$. The quantity $y_{ji}$ (the
number of letters $i$ that remain in $y_j$) is the minimum of
$x_{ji}$ (the number of letters $i$ that start out in $x_j$) and
the number of letters in $y_{j+1}$ that are available to block
these letters $i$ from moving to the other row. The latter
quantity is the total number of letters in $y_{j+1}$ that could
block letters $i$ (which is $\sum_{\alpha>i} y_{j+1,\alpha}$, the
number of letters in $y_{j+1}$ of value strictly greater than
$i$), minus the number of letters in $y_{j+1}$ that are already
being used to block letters in $x_j$ of value greater than $i$.
This last quantity is equal to the number of letters in $x_j$ that
are greater than $i$ and not moving to the other row, which is
$\sum_{\alpha>i} y_{j,\alpha}$.

This given, it is immediate that
\begin{equation} \label{eq:z}
\begin{split}
  z_{j+1,i} &= y_{j+1,i} + x_{ji} - y_{ji} \qquad \text{for $1\le
  j\le k-1$} \\
  z_{1,i} &= y_{1,i}.
\end{split}
\end{equation}

Continuing the computation of $\et{0}b$, one examines the first
row $y_1=z_1$. By the promotion method, one has
\begin{equation} \label{eq:prom eps0}
  \veps_0(b)=y_{11}.
\end{equation}
In particular,
\begin{equation} \label{eq:jeuzero}
  \et{0}b \not=0 \Leftrightarrow y_{11}>0.
\end{equation}
Suppose $\et{0}b\neq0$. As prescribed,
\begin{equation*}
  z'_{ji} = z_{ji} - \delta_{j,1}\delta_{i,1}.
\end{equation*}
Let $x'_{ji}=x(\et{0}b)$. Observe that one may use the same
process to pass from $\et{0}b$ to $z'$, as one does to pass from
$b$ to $z$. Therefore one may define $y'_j$ as in
\eqref{eq:twojeu} and $y'_{ji}$ as in \eqref{eq:y}. Then $z'$ is
related to $x'$ via $y'$ using \eqref{eq:z}.

We now prove by induction that there exists a unique sequence of
values $1=c_0<c_1<\dotsm<c_k=n+1$ such that
\begin{equation} \label{eq:y diff}
  y'_{ji} = y_{ji} - \delta_{i,c_{j-1}}.
\end{equation}
It is true by definition for $j=1$. By induction suppose
$c_0<\dotsc<c_{r-1}$ have been defined for some $r\ge 1$ and
\eqref{eq:y diff} holds for $1\le j \le r$. Write $y_r = u c_{r-1}
v$ for some words $u$ and $v$ such that $u$ does not contain
$c_{r-1}$. Then $y'_r=u v$. Consider \eqref{eq:twojeu} for $j=r$.
It can be interpreted that one row inserts the letters of $y_r$
from left to right, into the single-row tableau $z_{r+1}$, to
obtain the two-row tableau $y_{r+1} x_r$. Analogous statements
hold for the primed counterparts. By definition
$z'_{r+1}=z_{r+1}$. Clearly when inserting $y'_r=uv$ into
$z'_{r+1}$ the subword $u$ displaces the same letters as when
$y_r=u c_{r-1} v$ is inserted into $z_{r+1}=z'_{r+1}$. However the
insertion of $c_{r-1}$ is skipped, and one proceeds to insert $v$.
It follows from the definition of row insertion, that the multiset
of letters bumped by the insertion of the subword $v$ of $y'_r$,
is a submultiset of the letters bumped by the insertion of the
subword $v$ of $y_r$. The difference of these multisets is a
single letter $c_r$, which is greater than $c_{r-1}$ since it must
have been displaced by a letter in the subword $c_{r-1}v$ of
$y_r$. Thus \eqref{eq:y diff} holds for $j=r$, finishing the
induction.

Note that \eqref{eq:y diff} implies that
\begin{equation} \label{eq:x diff}
x'_{ji} = x_{ji} - \delta_{i,c_{j-1}} + \delta_{i,c_j}.
\end{equation}

Let $\chi(P)$ be $0$ or $1$ according as the statement $P$ is
false or true. Note that $y'$ is defined in terms of $x'$ as in
\eqref{eq:y}. Substituting \eqref{eq:x diff} and \eqref{eq:y diff}
one has
\begin{equation*}
\begin{split}
y_{ji}-\delta_{i,c_{j-1}}\!&=
\min(x_{ji}-\delta_{i,c_{j-1}}\!+\delta_{i,c_j}, \sum_{\alpha>i}
(y_{j+1,\alpha}-\delta_{\alpha,c_j})-(y_{j,\alpha}
-\delta_{\alpha,c_{j-1}}))
\\ &= \min(x_{ji}-\delta_{i,c_{j-1}}\!+\delta_{i,c_j},
\chi(c_{j-1}>i)-\chi(c_j>i)+\!\sum_{\alpha>i}
(y_{j+1,\alpha}-y_{j,\alpha}))
\end{split}
\end{equation*}
Adding $\delta_{i,c_{j-1}}$ we have
\begin{equation} \label{eq:new y}
y_{ji}=\min(x_{ji}+\delta_{i,c_j}, \chi(c_{j-1}\ge
i)-\chi(c_j>i)+\sum_{\alpha>i} (y_{j+1,\alpha}-y_{j,\alpha}))
\end{equation}
In the case $c_{j-1}<i<c_j$, \eqref{eq:new y} becomes
\begin{equation*}
y_{ji} = \min(x_{ji},-1+\sum_{\alpha>i} (y_{j+1,\alpha}-y_{j,\alpha})).
\end{equation*}
Comparing with \eqref{eq:y}, it follows that
\begin{equation}\label{eq:y between}
  y_{ji} = x_{ji}<\sum_{\alpha>i} (y_{j+1,\alpha}-y_{j,\alpha})
   \qquad\text{for $c_{j-1}< i<c_j$.}
\end{equation}
In the case $i=c_j$, \eqref{eq:new y} becomes
\begin{equation*}
  y_{j,c_j} =\min(x_{j,c_j}+1,\sum_{\alpha>i}
(y_{j+1,\alpha}-y_{j,\alpha})).
\end{equation*}
Comparing with \eqref{eq:y}, we have
\begin{equation} \label{eq:y end}
y_{j,c_j} = \sum_{\alpha>c_j} (y_{j+1,\alpha}-y_{j,\alpha}).
\end{equation}

In light of \eqref{eq:y between} and \eqref{eq:y end}, it follows
that the indices $c_j$ may be defined by the rule that
\begin{align}
c_0&=1 \\
\notag \text{for $j>1$, } & \text{ $c_j$ is the minimum index
$i>c_{j-1}$ such that}
\\
\label{eq:cj def}  y_{j,i} &= \sum_{\alpha>i}
(y_{j+1,\alpha}-y_{j,\alpha}).
\end{align}

\subsection{KKMMNN conditions imply promotion equations}
\label{ssec:KMNjeu} Let $b\in B^{k,l}$, $x=x(b)$, and $c\in C$
such that \eqref{eq:mas} and \eqref{eq:mass1} hold. Define $y$ in
terms of $x$ by \eqref{eq:y}. It is shown in this subsection that $y$
must satisfy \eqref{eq:y between} and \eqref{eq:y end}.

We will show by descending induction on $1\le j<k$ that
\begin{align}
\label{eq:x end}
x_{j,c_j} &\ge \sum_{\alpha> c_j} (y_{j+1,\alpha} - y_{j,\alpha}) &\\
\label{eq:x between}
x_{j,m} &< \sum_{\alpha>m} (y_{j+1,\alpha} - y_{j,\alpha})
 & \text{for $c_{j-1}<m<c_j$}.
\end{align}
By definition \eqref{eq:y} this implies immediately that
\begin{align}
\label{eq:y end'}
y_{j,c_j} &= \sum_{\alpha>c_j} (y_{j+1,\alpha} - y_{j,\alpha}) &\\
\label{eq:y between'}
y_{j,m} &= x_{j,m} & \text{for $c_{j-1}<m<c_j$,}
\end{align}
which also holds for $j=k$ by definition.

Let us prove \eqref{eq:x end} and \eqref{eq:x between} for
$j$ assuming that \eqref{eq:x end}-\eqref{eq:y between'} hold
for all labels greater than $j$. Define $s$ and $r_p$ for
$j\le p\le s$ as follows. Let $r_p>r_{p-1}$ ($r_j>c_j$ for $p=j$) be
minimal such that
\begin{equation}\label{eq:y rp}
y_{p,r_p}=\sum_{\alpha>r_p} (y_{p+1,\alpha}-y_{p,\alpha}),
\end{equation}
$r_p\ge c_{p+1}$ and $r_s<c_{s+1}$. The minimality of $r_p$ implies
\begin{equation}\label{eq:y rp between}
y_{p,m}=x_{p,m} \qquad \text{for $r_{p-1}<m<r_p$ ($c_j<m<r_j$ for $p=j$)}.
\end{equation}

First we show by induction on $p=j,j+1,\ldots,s$ that \eqref{eq:x end}
is equivalent to
\begin{equation}\label{eq:inter}
\sum_{\beta=j+1}^p \; \sum_{c_{\beta-1}<\alpha<c_{\beta}} x_{\beta,\alpha}
 + \sum_{c_p<\alpha\le r_p} y_{p+1,\alpha}
\le \sum_{c_j\le \alpha<r_j} x_{j,\alpha} + \sum_{\beta=j+1}^p\;
 \sum_{r_{\beta-1}<\alpha<r_\beta} x_{\beta,\alpha}.
\end{equation}
The case $p=j$ follows immediately from \eqref{eq:x end} by
using \eqref{eq:y rp} and \eqref{eq:y rp between}.
Now assume that \eqref{eq:inter} holds for $p-1\ge j$. The left hand side
of \eqref{eq:inter} at $p-1$ can then be transformed as follows:
\begin{equation*}
\begin{split}
l.h.s.&=
 \sum_{\beta=j+1}^p \; \sum_{c_{\beta-1}<\alpha<c_{\beta}} x_{\beta,\alpha}
 + \sum_{c_p\le \alpha\le r_{p-1}} y_{p,\alpha}\\
&= \sum_{\beta=j+1}^p \; \sum_{c_{\beta-1}<\alpha<c_{\beta}} x_{\beta,\alpha}
 + \sum_{\alpha>c_p} y_{p+1,\alpha} - \sum_{\alpha>r_{p-1}} y_{p,\alpha}\\
&= \sum_{\beta=j+1}^p \; \sum_{c_{\beta-1}<\alpha<c_{\beta}} x_{\beta,\alpha}
 - \sum_{r_{p-1}<\alpha<r_p} x_{p,\alpha}
 + \sum_{c_p<\alpha\le r_p} y_{p+1,\alpha},
\end{split}
\end{equation*}
where the first equality follows from \eqref{eq:y between'} at $p$,
for the second equality follows from \eqref{eq:y end'} at $p$,
and for the  last equality we used \eqref{eq:y rp} and \eqref{eq:y rp between}.
Hence \eqref{eq:inter} holds at $p$.

For $p=s$, all $y_{s+1,\alpha}$ in \eqref{eq:inter} can be replaced by
$x_{s+1,\alpha}$ by \eqref{eq:y between'} since by assumption $r_s<c_{s+1}$.
Hence \eqref{eq:inter} with $p=s$ is equivalent to
\begin{equation*}
\Delta(c)\le \Delta(c_0,c_1,\ldots,c_{j-1},r_j,\ldots,r_s,c_{s+1},\ldots,c_k)
\end{equation*}
which is true by \eqref{eq:mas}. This proves \eqref{eq:x end}.

Next we prove \eqref{eq:x between} by induction on $m$. By induction
hypothesis we may assume that \eqref{eq:x between} holds for all indices
greater than $m$.
Define $s$ and $r_p$ for $j\le p\le s$ as follows.
Set $r_j=m$ and let $r_{p-1}<r_p<c_p$ be minimal such that
\begin{equation} \label{eq:y rp2}
y_{p,r_p}=\sum_{\alpha>r_p} (y_{p+1,\alpha} - y_{p,\alpha})
\end{equation}
holds. Let $s<k-1$ be minimal such that $r_{s+1}$ does not exist;
otherwise let $s=k-1$. The minimality of $r_p$ implies
\begin{equation*}
y_{p,m}=x_{p,m} \qquad \text{for $r_{p-1}<m<r_p$}.
\end{equation*}
We prove by induction on $p$ that
\eqref{eq:x between} is equivalent to
\begin{equation}\label{eq:inter1}
\sum_{\beta=j}^{p+1} \; \sum_{c_{\beta-1}<\alpha<c_\beta} x_{\beta,\alpha}
< \sum_{c_{j-1}<\alpha<r_j} x_{j,\alpha}
 +\sum_{\beta=j+1}^p \; \sum_{r_{\beta-1}<\alpha<r_\beta} x_{\beta,\alpha}
 +\sum_{r_p<\alpha<c_{p+1}} y_{p+1,\alpha}.
\end{equation}
To see that this equivalence holds at $p=j$, we use the induction
hypothesis $y_{j,\alpha}=x_{j,\alpha}$ for $m<\alpha<c_j$ and
\eqref{eq:y end'}. With this \eqref{eq:x between} reads
\begin{equation*}
\sum_{r_j\le \alpha< c_j} x_{j,\alpha}
 < \sum_{r_j<\alpha\le c_j} y_{j+1,\alpha}.
\end{equation*}
Adding
\[
\sum_{c_{j-1}<\alpha<r_j} x_{j,\alpha} +\sum_{c_j<\alpha<c_{j+1}}
x_{j+1,\alpha} =\sum_{c_{j-1}<\alpha<r_j}
x_{j,\alpha}+\sum_{c_j<\alpha<c_{j+1}}
 y_{j+1,\alpha}
\]
to both sides yields \eqref{eq:inter1} at $p=j$. Now assume that
\eqref{eq:inter1} holds for $p-1\ge j$. The last term in this
inequality can be rewritten as
\begin{equation*}
\begin{split}
\sum_{r_{p-1}<\alpha<c_p} y_{p,\alpha} &=
 \sum_{r_{p-1}<\alpha<r_p} y_{p,\alpha}
 +\sum_{\alpha>r_p}(y_{p+1,\alpha}-y_{p,\alpha})
 +\sum_{r_p<\alpha<c_p} y_{p,\alpha}\\
&=\sum_{r_{p-1}<\alpha<r_p} x_{p,\alpha}
 +\sum_{\alpha>r_p} y_{p+1,\alpha}
 -\sum_{\alpha\ge c_p} y_{p,\alpha}\\
&=\sum_{r_{p-1}<\alpha<r_p} x_{p,\alpha}
 +\sum_{r_p<\alpha\le c_p} y_{p+1,\alpha}.
\end{split}
\end{equation*}
Using this and adding $\sum_{c_p<\alpha<c_{p+1}} x_{p+1,\alpha}
=\sum_{c_p<\alpha<c_{p+1}} y_{p+1,\alpha}$ to both sides of \eqref{eq:inter1}
at $p-1$ yields \eqref{eq:inter1} at $p$.

At $p=s$ the last term in \eqref{eq:inter1} can be replaced by
$\sum_{r_s<\alpha<c_{s+1}} x_{s+1,\alpha}$. For if $s<k-1$ then
there is no $r_{s+1}<c_{s+1}$ such that \eqref{eq:y rp2} holds, and
if $s=k-1$ then by \eqref{eq:y}, $x_{k,\alpha}=y_{k,\alpha}$ for
all $1\le\alpha\le n$. Hence \eqref{eq:inter1} at $p=s$ is
equivalent to
\begin{equation*}
\Delta(c)<\Delta(c_0,\ldots,c_{j-1},r_j,\ldots,r_s,c_{s+1},\ldots,c_k)
\end{equation*}
which is true by \eqref{eq:mass1} since $c_{j-1}<r_j<c_j$.

\subsection{Equivalence of the two methods for $\et{0}$}
\label{ssec:equivzero} Fix $b\in B^{k,\ell}$ and let $x=x(b)$. Let
$c$ be defined as in \eqref{eq:mas} and \eqref{eq:mass1}. Let $y$
be defined by \eqref{eq:y}. By subsection \ref{ssec:KMNjeu}, $y$
satisfies \eqref{eq:y between} and \eqref{eq:y end}. We show that
the KKMMNN and promotion methods of computing $\et{0}b$, are
equivalent.

First, we show the equivalence of the conditions given in
subsections \ref{ssec:KMN} and \ref{ssec:jeu}, that
$\et{0}b\neq0$. Explicitly, we show that $c\not=c_{\min}$ or
$x_{k,k}>0$, if and only if $y_{11}>0$.

Recall that $y_1$ has length $\ell-x_{k,n+1}=\sum_{i\ge1} y_{1i}$.
Applying \eqref{eq:y between} and \eqref{eq:y end} repeatedly, we
have
\begin{equation} \label{eq:y11}
\begin{split}
  \ell-x_{k,n+1}-y_{11} &= \sum_{i>1} y_{1i} \\
  &= \sum_{c_0<i<c_1} x_{1i} + \sum_{i\ge c_1} y_{1i} \\
  &= \sum_{c_0<i<c_1} x_{1i} + \sum_{i>c_1} y_{2i} \\
  &= \sum_{c_0<i<c_1} x_{1i} + \sum_{c_1<i<c_2} x_{2i} +
  \sum_{i\ge c_2} y_{2i} \\
  &= \dotsm = \Delta(c).
\end{split}
\end{equation}

By the definition of $c_{\min}$, $\Delta(c_{\min}) =
\sum_{i=k+1}^n x_{ki}$. The $k$-th row of $x$ has length $\ell$.
By semistandardness it can only have elements in the set
$\{k,k+1,\dotsc,n+1\}$. Therefore
$\Delta(c_{\min})=\ell-x_{kk}-x_{k,n+1}$. By \eqref{eq:y11},
\begin{equation}
  y_{11} =
  \ell-x_{k,n+1}-\Delta(c)=x_{kk}+\Delta(c_{\min})-\Delta(c).
\end{equation}
In light of \eqref{eq:mass1} the desired equivalence is evident.

We are now left to show that both methods agree provided that
$\et{0}b\neq0$. Assume $\et{0}b\neq0$. Let $c'\in C$ be the
sequence defined by the promotion method. Then $y$ (which is
defined only in terms of $x$) satisfies \eqref{eq:y between} and
\eqref{eq:y end} for both $c$ and $c'$. It follows that $c=c'$ so
that the two methods agree.

\subsection{Proof of \eqref{eq:Bkl eps0} and \eqref{eq:Bkl phi0}}
Equation \eqref{eq:Bkl eps0} follows immediately from
\eqref{eq:y11}, \eqref{eq:prom eps0}, and the equivalence of the
KKMMNN and promotion methods for $\et{0}$.

Equation \eqref{eq:i-wt} implies
\begin{equation}
  \vphi_0(b)-\veps_0(b)=\langle h_0, \wt\, b \rangle
  =x_{k,n+1}-x_{11}.
\end{equation}
Equation \eqref{eq:Bkl phi0} follows from this and \eqref{eq:Bkl
eps0}. Note that although different sequences $c\in C^x_{\min}$
are used in the definition of $\et{0}$ and $\ft{0}$, they both
attain the same minimum value of the function $\Delta$.

\subsection{Equivalence of the two methods for $\ft{0}$}
\label{ssec:ft0} The promotion definition of $\ft{0}$ is
essentially the inverse of that of $\et{0}$. However the KKMMNN
definition of $\ft{0}$ is not defined so as to be obviously equal
to the inverse of $\et{0}$.

Let $b'\in B^{k,\ell}$ be such that $\ft{0}b'=b\not=0$. Then of
course $b'=\et{0}b$. Write $x'=x(b')$ and $x=x(b)$. Let $c$ be the
unique $\subseteq$-minimum element of $C^x_{\min}$. Then $x'$ and
$x$ are related as in \eqref{eq:Bkl e0}.

It must be shown that the KKMMNN method for $\ft{0}$ sends $b'$ to
$b$. Let $c'$ and $c''$ be the unique $\subseteq$-maximum and
minimum elements of $C^{x'}_{\min}$ respectively; in particular
$\Delta_{x'}(c')=\Delta_{x'}(c'')$. By \eqref{eq:Bkl f0} it must
be shown that $c'=c$.

To prove (1) for $\ft{0} b'$, note that $\ell$ is the length of
the first row of $b'$ and that $x'_{1i}=0$ for $i>k'+1$ by
semistandardness. The definition of $c_{\max}$ gives
\begin{equation*}
  \ell = \sum_{i=1}^{k'+1} x'_{1i} = x'_{11}+\Delta_{x'}(c_{\max})+x'_{1,k'+1}.
\end{equation*}
By \eqref{eq:Bkl phi0} we have
\begin{equation*}
  \vphi_0(b') = x'_{1,k'+1}+\Delta_{x'}(c_{\max})  - \Delta_{x'}(c').
\end{equation*}
{}From this and the definition of $c'\in C^{x'}_{\min}$ it is seen that
the conditions for $\ft{0}b'=0$ given by the two methods
coincide.

By \eqref{eq:Bkl eps0} applied to both $b$ and $b'$ we have
\begin{equation} \label{eq:deltas}
\begin{split}
 \Delta_x(c) &= -\veps_0(b)+\ell-x_{k,n+1}
 = -(\veps_0(b')+1)+\ell-(x'_{k,n+1}-1) \\
 &= \Delta_{x'}(c'') = \Delta_{x'}(c').
\end{split}
\end{equation}

By \eqref{eq:Bkl e0} it follows that for all $m\in C$
\begin{equation} \label{eq:newdelta}
\Delta_{x'}(m)=\Delta_x(m) + \sum_{j=1}^k
(\chi(m_{j-1}<c_j<m_j)-\chi(m_{j-1}<c_{j-1}<m_j)).
\end{equation}
Evaluating \eqref{eq:newdelta} at $m=c$ and applying
\eqref{eq:deltas} we have $\Delta_{x'}(c)=\Delta_{x'}(c')$.
Evaluating \eqref{eq:newdelta} at $m\supseteq c$ such that
$m\not=c$, we see that $\Delta_{x'}(m) \ge \Delta_x(m) >
\Delta_x(c)$ by \eqref{eq:mass1}. In other words, $c$ is the
unique $\subseteq$-maximum element of $C^{x'}_{\min}$, that is,
$c=c'$.

\end{document}